\documentclass[12pt,a4paper]{amsart}
\usepackage{mathtools}
\usepackage{xcolor}
\usepackage{amsmath,amstext,amssymb,amsthm,amsfonts}
\usepackage{graphicx}
\usepackage[all]{xy} \xyoption{poly}
\usepackage{ mathrsfs }

\newcommand{\nc}{\newcommand}
\nc{\one}{\mbox{\bf 1}}
\nc{\invtensor}{\underset{\leftarrow}{\otimes}}
\nc{\const}{\operatorname{const}}

\nc{\ad}{\operatorname{ad}}
\nc{\ev}{\operatorname{ev}}
\nc{\tr}{\operatorname{tr}}
\nc{	\Gr}{\mathscr{K}}
\nc{\rGr}{\operatorname{rGr}}
\nc{\atyp}{\operatorname{atyp}}
\nc{\tp}{\operatorname{top}}
\nc{\rank}{\operatorname{rank}}
\nc{\corank}{\operatorname{corank}}
\nc{\codim}{\operatorname{codim}}
\nc{\sdim}{\operatorname{sdim}}
\nc{\mult}{\operatorname{mult}}
\nc{\ds}{\operatorname{ds}}
\nc{\tail}{\operatorname{tail}}
\nc{\howl}{\operatorname{howl}}
\nc{\spn}{\operatorname{span}}
\nc{\Sym}{\operatorname{Sym}}
\nc{\Core}{\operatorname{Core}}
\nc{\id}{\operatorname{id}}
\nc{\Id}{\operatorname{Id}}
\nc{\Ree}{\operatorname{Re}}
\nc{\hi}{\operatorname{hi}}
\nc{\htt}{\operatorname{ht}}
\nc{\at}{\operatorname{at}}
\nc{\str}{\operatorname{str}}
\nc{\Iso}{\operatorname{Iso}}
\nc{\Ker}{\operatorname{Ker}}
\nc{\rker}{\operatorname{rKer}}
\nc{\im}{\operatorname{Im}}
\nc{\osp}{\mathfrak{osp}}
\nc{\sgn}{\operatorname{sgn}}
\nc{\F}{\operatorname{F}}
\nc{\Mod}{\operatorname{Mod}}
\nc{\DS}{\operatorname{DS}}
\nc{\Soc}{\operatorname{Soc}}
\nc{\Inj}{\operatorname{Inj}}
\nc{\Hom}{\operatorname{Hom}}
\nc{\End}{\operatorname{End}}
\nc{\supp}{\operatorname{supp}}
\nc{\smult}{\operatorname{smult}}
\nc{\sOmega}{\operatorname{s}\Omega}
\nc{\Card}{\operatorname{Card}}
\nc{\Ann}{\operatorname{Ann}}
\nc{\Arc}{\operatorname{Arc}}
\nc{\arc}{\operatorname{arc}}
\nc{\Ind}{\operatorname{Ind}}
\nc{\Coind}{\operatorname{Coind}}
\nc{\wt}{\operatorname{hwt}}
\nc{\hwt}{\operatorname{wt}}
\nc{\rk}{\operatorname{rank}}
\nc{\ch}{\operatorname{ch}}
\nc{\sch}{\operatorname{sch}}
\nc{\mdim}{\operatorname{mdim}}
\nc{\Stab}{\operatorname{Stab}}

\nc{\Irr}{\operatorname{Irr}}
\nc{\Spec}{\operatorname{Spec}}
\nc{\Res}{\operatorname{Res}}
\nc{\res}{\operatorname{res}}
\nc{\Aut}{\operatorname{Aut}}
\nc{\Ext}{\operatorname{Ext}}
\nc{\Prec}{\operatorname{Prec}}
\nc{\Fract}{\operatorname{Fract}}
\nc{\gr}{\operatorname{gr}}
\nc{\diag}{\operatorname{diag}}
\nc{\deff}{\operatorname{def}}
\nc{\core}{\operatorname{core}}
\nc{\HC}{\operatorname{HC}}
\nc{\Ch}{\operatorname{Ch}}
\nc{\dpth}{\operatorname{dpth}}
\nc{\sw}{\operatorname{sw}}
\nc{\red}{\operatorname{red}}
\nc{\pari}{\operatorname{par}}
\nc{\pos}{\operatorname{pos}}
\nc{\Cl}{\mathcal{C}\ell}

\nc{\wdchi}{\widetilde{\chi}}
\nc{\wdH}{\widetilde{H}}
\nc{\wdN}{\widetilde{N}}
\nc{\wdM}{\widetilde{M}}
\nc{\wdO}{\widetilde{O}}
\nc{\wdR}{\widetilde{R}}

\nc{\wdV}{\widetilde{V}}

\nc{\wdC}{\widetilde{C}}

\nc{\zero}{\operatorname{zero}}
\nc{\nonzero}{\operatorname{nonzero}}

\nc{\Nonzero}{\operatorname{Nonzero}}

\nc{\Obj}{\operatorname{Obj}}
\nc{\Dglie}{\operatorname{{\mathcal D}glie}}
\nc{\Fin}{\operatorname{{\mathcal F}in}}
\nc{\pr}{\operatorname{pr}}
\nc{\Adm}{\operatorname{\mathcal{A}dm}}
\nc{\fg}{\mathfrak g}

\nc{\Sg}{{\cS(\fg)}}
\nc{\Shg}{{\cS(\fhg)}}
\nc{\Ug}{{\cU(\fg)}}
\nc{\Uhg}{{\cU(\fhg)}}
\nc{\Sh}{{\cS(\fh)}}
\nc{\Uh}{{\cU(\fh)}}
\nc{\Uhh}{{\cU(\fhh)}}
\nc{\Zg}{{{\mathcal{Z}}(\fg)}}

\nc{\Vir}{{\mathcal{V}ir}}
\nc{\NS}{{\mathcal{N}S}}

\nc{\tZg}{{\widetilde{\mathcal Z}({\mathfrak g})}}
\nc{\Zk}{{\mathcal Z}({\mathfrak k})}

\newcommand{\ZZ}{\mathbb{Z}}

\nc{\Up}{{\mathcal U}({\mathfrak p})}
\nc{\Ah}{{\mathcal A}({\mathfrak h})}
\nc{\Ag}{{\mathcal A}({\mathfrak g})}
\nc{\Ap}{{\mathcal A}({\mathfrak p})}
\nc{\Zp}{{\mathcal Z}({\mathfrak p})}
\nc{\cR}{\mathcal R}
\nc{\cS}{\mathcal S}
\nc{\cP}{\mathcal P}
\nc{\cT}{\mathcal{T}}
\nc{\CC}{\mathcal C}
\nc{\cA}{\mathcal A}
\nc{\cV}{\mathcal V}
\nc{\cU}{\mathcal U}
\nc{\cZ}{\mathcal Z}
\nc{\cM}{\mathcal M}
\nc{\cL}{\mathcal L}
\nc{\cF}{\mathcal F}

\nc{\cB}{\mathcal{B}}

\nc{\fo}{\mathfrak o}

\nc{\fa}{\mathfrak a}

\nc{\CO}{\mathcal O}
\nc{\CR}{\mathcal R}

\nc{\cK}{\mathcal{K}}
\nc{\cW}{\mathcal{W}}
\nc{\bM}{\mathbf{M}}
\nc{\bL}{\mathbf{L}}
\nc{\bN}{\mathbf{N}}

\nc{\zq}{\mathpzc q}

\nc{\fl}{\mathfrak l}
\nc{\fn}{\mathfrak n}
\nc{\fm}{\mathfrak m}
\nc{\fp}{\mathfrak p}
\nc{\fh}{\mathfrak h}
\nc{\ft}{\mathfrak t}
\nc{\fk}{\mathfrak k}
\nc{\fb}{\mathfrak b}
\nc{\fs}{\mathfrak s}
\nc{\fB}{\mathfrak B}

\nc{\vareps}{\varepsilon}
\nc{\varesp}{\varepsilon}
\nc{\veps}{\varepsilon}

\nc{\fsl}{\mathfrak{sl}}
\nc{\fgl}{\mathfrak{gl}}
\nc{\fso}{\mathfrak{so}}
\nc{\fosp}{\mathfrak{osp}}
\nc{\fsp}{\mathfrak{sp}}
\nc{\fq}{\mathfrak q}
\nc{\fsq}{\mathfrak{sq}}
\nc{\fpsq}{\mathfrak{psq}}
\nc{\fpq}{\mathfrak{pq}}

\def\d{\partial}
\def\sub{\subseteq}

\nc{\fhg}{\hat{\fg}}
\nc{\fhn}{\hat{\fn}}
\nc{\fhh}{\hat{\fh}}
\nc{\fhb}{\hat{\fb}}
\nc{\hrho}{\hat{\rho}}

\nc{\hsl}{\hat{\fsl}}
\nc{\fpo}{\mathfrak{po}}
\nc{\dirlim}{\underset{\rightarrow}{\lim}\,}
\nc{\nen}{\newenvironment}
\nc{\ol}{\overline}
\nc{\ul}{\underline}
\nc{\ra}{\rightarrow}
\nc{\lra}{\longrightarrow}
\nc{\Lra}{\Longrightarrow}
\nc{\bo}{\bar{1}}
\nc{\Lla}{\Longleftarrow}

\nc{\Llra}{\Longleftrightarrow}

\nc{\thla}{\twoheadleftarrow}

\nc{\lang}{(}
\nc{\rang}{)}

\nc{\hra}{\hookrightarrow}

\nc{\iso}{\overset{\sim}{\lra}}

\nc{\ssubset}{\underset{\not=}{\subset}}

\nc{\vac}{|0\rangle}

\nc{\simka}{{\ \scriptscriptstyle _{{\sim}}^\text{\tiny{k}}\ }}

\nc{\Thm}[1]{Theorem~\ref{#1}}
\nc{\Prop}[1]{Proposition~\ref{#1}}
\nc{\Lem}[1]{Lemma~\ref{#1}}
\nc{\Cor}[1]{Corollary~\ref{#1}}
\nc{\Conj}[1]{Conjecture~\ref{#1}}
\nc{\Claim}[1]{Claim~\ref{#1}}
\nc{\Defn}[1]{Definition~\ref{#1}}
\nc{\Exa}[1]{Example~\ref{#1}}
\nc{\Rem}[1]{Remark~\ref{#1}}
\nc{\Note}[1]{Note~\ref{#1}}
\nc{\Quest}[1]{Question~\ref{#1}}
\nc{\Hyp}[1]{Hypoth\`ese~\ref{#1}}
\nen{thm}[1]{\label{#1}{\bf Theorem.\ } \em}{}
\nen{prop}[1]{\label{#1}{\bf Proposition.\ } \em}{}
\nen{lem}[1]{\label{#1}{\bf Lemma.\ } \em}{}
\nen{cor}[1]{\label{#1}{\bf Corollary.\ } \em}{}
\nen{conj}[1]{\label{#1}{\bf Conjecture.\ } \em}{}

{\newtheorem{theorem}{Theorem} 
	\newtheorem{lemma}{Lemma}
	\newtheorem{proposition}{Proposition}
	\newtheorem{corollary}{Corollary}
	
	\theoremstyle{definition}
	
	\theoremstyle{definition}
	\newtheorem{remark}{Remark}
	\theoremstyle{definition}
	}
\numberwithin{theorem}{section}
\numberwithin{lemma}{section}
\numberwithin{remark}{section}
\numberwithin{definition}{section}
\numberwithin{corollary}{section}
\numberwithin{proposition}{section}
\numberwithin{example}{section}

\begin{document}
\setcounter{section}{0}
\setcounter{tocdepth}{1}

\title[Grothendieck rings]{On the Grothendieck ring of a quasireductive Lie superalgebra}
\author{ M. Gorelik, V. Serganova, A. Sherman}

\address{Weizmann Institute of Science, Rehovot, Israel}
\email{maria.gorelik@gmail.com}

\address{Dept.~of Mathematics, University of California, Berkeley, USA}
\email{serganov@math.berkeley.edu}
	
\address{School of Mathematics and Statistics, University of Sydney, Sydney, Australia} 
\email{xandersherm@gmail.com}

\date{}

\begin{abstract} Given a Lie superalgebra $\fg$ and a maximal quasitoral subalgebra $\fh$, we consider properties of restrictions of $\fg$-modules to $\fh$. This is a natural generalization of the study of characters in the case when $\fh$ is an even maximal torus.  We study the case of $\fg=\fq_n$ with $\fh$ a Cartan subalgebra, and prove several special properties of the restriction in this case, including an explicit realization of the $\fh$-supercharacter ring.
\end{abstract}


\medskip


\maketitle

\section{Introduction}

\subsection{Maximal toral subalgebras and restriction}\label{section-maxtoral}
 Let $\fg$ be a Lie superalgebra, not necessarily finite-dimensional.  Assume that
$\fg$ contains a subalgebra $\fh$ such that

\begin{enumerate}
	\item $\ft:=\fh_{\ol{0}}$ has diagonalizable adjoint action on $\fg$ with finite-dimensional weight spaces; and,
	\item  $\fh=\fg^{\ft}$, where $\fg^{\ft}$ denotes the centralizer of $\ft$ in $\fg$;
\end{enumerate} 
We will call subalgebras with the above properties \emph{maximal quasitoral}. Maximal quasitoral subalgebras play an analogous role to maximal toral subalgebras of Lie algebras.  In the purely even setting, the restriction of a representation to a maximal toral subalgebra is exactly the data of its character.  The character of a representation is a powerful invariant, and providing (nice) formulas for characters of irreducible representations is a central problem in representation theory.

For many Lie superalgebras of interest, maximal quasitoral subalgebras are even, i.e. $\fh=\ft$ (e.g.~for $\fg\fl_{m|n},\fo\fs\fp_{m|2n},\fp_n,\dots$ See \cite{KLie} for the definition of these Lie superalgebras).  In this case the restriction of an irreducible highest weight representation to $\ft$ completely determines it.  But when $\fh\neq\ft$, as in the case of the queer Lie superalgebra $\fq_n$, this is no longer holds, as for certain simple modules $L$ one has
$\Res^{\fg}_{\ft} L\cong \Res^{\fg}_{\ft}\Pi L$ and so the restriction map of Grothendieck groups $\Gr(\fg)\to\Gr(\ft)$ is not injective. 

On the other hand the restriction to $\fh$ is fine enough to make such distinctions.  Choose a Borel subalgebra $\fb$ of $\fg$ containing $\fh$ (from a triangular decomposition as explained in Section \ref{section-maxl-qtoral-subalg}), and let $\CC$ be the full subcategory of $\fg$-modules such that each module in $\CC$ is of 
finite length, has a diagonal action of $\ft$, and is locally finite over $\fb$ (note these conditions can be slightly weakened, see Corollary \ref{corollary-res-Gr-embedding}).  Write $\cM^s(\fh)$ for the category of $\fh$-modules with diagonal $\ft$ action and finite-dimensional weight spaces.  Then the restriction functor $\Res^{\fg}_{\fh}$
induces an injective map $\Gr(\CC)\to \Gr(\cM^s(\fh))$, see Corollary \ref{corollary-res-Gr-embedding}.

It is therefore of interest to understand the $\fh$-character of irreducible $\fg$-modules.  A difficulty that arises is that the Grothendieck ring of finite-dimensional $\fh$-modules does not have a simple manifestation like in the case of $\ft$-modules.  Nevertheless a description can be given, and its structure is interesting in its own right.

\subsection{$\Gr_-(\fg)$ and $\Gr_+(\fg)$}  Let $\Gr(\fg)$ denote the Grothendieck group of finite-dimensional $\fg$-modules, and write $M_{gr}$ for the natural image of a finite-dimensional module $M$ in $\Gr(\fg)$.  This ring admits two natural quotients: $\Gr_+(\fg)$, which is obtained by identifying $M_{gr}=(\Pi M)_{gr}$, and $\Gr_-(\fg)$, obtained by identifying $M_{gr}=-(\Pi M)_{gr}$.  There is a natural embedding of $\Gr(\fg)$ into $\Gr_+(\fg)\times\Gr_-(\fg)$, and this embedding becomes an isomorphism over $\mathbb{Q}$; thus a proper understanding of $\Gr_+(\fg)$ and $\Gr_-(\fg)$ suffices for the understanding of $\Gr(\fg)$ (see Section 5 for a precise relationship between $\Gr(\fg)$ and $\Gr_+(\fg)\times\Gr_-(\fg)$).

The ring $\Gr_+(\fg)$ behaves like a character ring, and in fact embeds into $\Gr_+(\ft)$ under the restriction map; thus the information it carries is less interesting from our standpoint.  On the other hand, $\Gr_-(\fg)$, the \emph{reduced Grothendieck ring}, has a nontrivial, even exotic structure as a ring, and will be our main object of study.  The restriction map $\Gr_-(\fg)\to\Gr_-(\ft)$ is sometimes very far from being an embedding; for instance if $\fg=\fq_n$ the image of any nontrivial finite-dimensional irreducible module is zero by \cite{Cheng}.  Thus is it necessary to study instead $\Gr_-(\fg)\to\Gr_-(\fh)$.

A further advantage of using $\Gr_-(\fg)$ is that the Duflo-Serganova functor, while not being exact, always induces a morphism $ds_x:\Gr_-(\fg)\to\Gr_-(\fg_x)$ for appropriate $x\in\fg_{\ol{1}}$.  It has been known for some time (see \cite{HR}) that for Kac-Moody superalgebras, the map induced by $ds_x$ on supercharacters is given by restriction to $\mathfrak{t}_x$, the Cartan subalgebra of $\fg_x$.  This is a reflection of a more general property of $ds_x$, discussed in Section 8, which shows that $ds_x$ can always be thought of as a restriction map to $\fg_x$.   Therefore in our setting the induced map $ds_x:\Gr_-(\fh)\to\Gr_-(\fh_x)$ is also given by restriction of modules from $\fh$ to $\fh_x$, where $\fh_x$ is a maximal quasitoral subalgebra of $\fg_x$. 

\subsection{Results for $\fg=\fq_n$}  Let $\fg=\fq_n$; for two weights $\lambda,\mu,$ write $\lambda\sim\mu$ if $L(\lambda)$ and $L(\mu)$ lie in the same block of $\Fin(\fq_n)$.  Then for a $\fg$-module $M$, write $\sch_{\fh}M$ for the natural image of $M$ in $\Gr_-(\Fin\fh)$. 

\begin{theorem}\label{theorem-intro_1}
	\begin{equation}
	\sch_{\fh} L(\lambda)=\sum_{\mu\sim\lambda} \sch_{\fh} L(\lambda)_{\mu}
	\end{equation}
	Equivalently, if $\mu\not\sim\lambda$ then $[L(\lambda):C(\mu)]=[L(\lambda):\Pi C(\mu)]$, where $C(\mu)$ is an irreducible $\fh$-module of weight $\mu$.
\end{theorem}

See Theorem \ref{theorem-sch-q-long} for a proof of this statement, and a stronger result.  In Section 6 we present a version of the above result which holds for a more general Lie superalgebra, but only for particular blocks, those which are the closest to being `typical'.

\subsubsection{Supercharacter isomorphism}  Let $\CC$ denote the subcategory of $\fh$-modules with weights which appear in some finite-dimensional representation of $\mathfrak{gl}_n=(\fq_n)_{\ol{0}}$.  The Weyl group $W=S_n$ has a natural action on $\Gr_-(\CC	)$ in this case, and we may consider the invariant subalgebra.  We prove that $\Gr_-(\CC)^W$ has a basis given by $\{a_{\lambda}\}_{\lambda\in P^+(\fg)}$, where $P^+(\fg)$ are the dominant weights for $\fg=\fq_n$, and
\[
a_{\lambda}=\sum\limits_{w\in W/\Stab_W(\lambda)}[C(\lambda)^w],
\]
where $(-)^w$ denotes the twisting functor.

It is clear that the supercharacter map induces an embedding $\sch_{\fh}:\Gr_-(\fg)\to\Gr_-(\CC)^W$; such a result holds for any quasireductive Lie superalgebra (see Section 7).  For $\fg=\fq_n$ we have (see Corollary \ref{corollary-red-Gr-q(n)}):

\begin{theorem}
		The map $\sch_{\fh}:\Gr_-(\fg)_{\mathbb{Q}}\to\Gr_-(\CC)_{\mathbb{Q}}^W$ is an isomorphism of rings.
\end{theorem}

In the above theorem, the subscript $\mathbb{Q}$ means that we extend scalars to $\mathbb{Q}$.  To obtain an isomorphism we only need to invert $2$, in fact.  This is a consequence of the work in Section 9.

\begin{remark}
The ring $\Gr_-(\fg)$ has a natural basis given by irreducible modules, and an important question is to understand the relation between this basis and the basis $\{a_{\lambda}\}$.
\end{remark}

\subsubsection{Realization of $\Gr_-(\cF(\fq_n)_{int})$} Using the above isomorphism, we are able to provide another realization of $\Gr_-(\cF(\fg)_{int})$, where $\cF(\fg)_{int}$ denotes the category of finite-dimensional $\fq_{n}$-modules with integral weights.  To simplify the explanation for the introduction, we will explain this realization over $\mathbb{C}$.  

Let $V:=\mathbb{C}^{\mathbb{Z}\setminus\{0\}}$ denote the complex vector space with basis $\{v_i\}_{i\in\mathbb{Z}\setminus\{0\}}$.  Then write
\[
A=\bigwedge V:=\bigoplus\limits_{n\in\mathbb{N}}{\bigwedge}^nV.
\]
This is a superalgebra, and we write $A_{\ol{0}}=\bigoplus\limits_{n\in\mathbb{N}}\bigwedge^{2n}V$ for the even part. Let $J_k\sub A$ be the ideal generated by $\bigwedge^{k+1}V$.  For the following, see Theorem \ref{thm 9.3}.

\begin{theorem}
	We have an explicit isomorphism of algebras
	\[
	\Gr_-(\cF(\fq_n)_{int})\otimes_{\mathbb{Z}}\mathbb{C}\to (A/J_n)_{\ol{0}}.
	\]
	where, up to scalar, $a_{\lambda}$ is mapped to $v_{\lambda}:=v_{j_1}\wedge\cdots\wedge v_{j_k}$, where $j_1,\dots,j_k$ are the nonzero coordinates of $\lambda$.
\end{theorem}
 
\subsubsection{Relation to the Duflo-Serganova functor}  
For $\fq_n$, the maps $ds_x$ depend only on the \emph{rank} of $x$, which is a non-negative half-integer $s\in\frac{1}{2}\mathbb{N}$.  Thus we write $ds_x=ds_s$, where $s=\rank(x)$, and this gives a map $ds_s:\Gr_-(\cF(\fq_n))\to\Gr_-(\cF(\fq_{n-2s}))$.  We have the following simple formula for $ds_s$ in terms of the basis $\{a_{\lambda}\}$, using the previously noted fact that $ds_x$ is given by restriction:
\[
\ds_s(a_{\mu})=\left\{\begin{array}{lcl}
0 & & \text{ if }
\zero \mu<2s\\
a_{\mu'} & & \text{ if }
\zero \mu\geq 2s\end{array}\right.
\]
where $\zero\mu$ is the number of zero coordinates that $\mu$ has, and $\mu'\in P^+(\fg_x)$ is such that $\mu'$ and $\mu$ have the same nonzero coordinates.  In \cite{GS}, we compute $ds_s$ on the basis $\{[L(\lambda)]\}_{\lambda\in P^+(\fg)}$ of irreducible modules; remarkably it admits a similarly simple expression:
\[
\ds_s([L(\lambda)])=\left\{\begin{array}{lcl}
	0 & & \text{ if }
	\zero \lambda<2s, \\
{[}L(\lambda')] & & \text{ if }
	\zero \lambda\geq 2s.\end{array}\right.
\]

%\subsection{Summary of sections}  In Section 2 we introduce the class of Lie superalgebras $\fg$ we consider as well as introduce the Clifford superalgebra and its irreducible representations.  In Section 3 we describe the representation theory of quasitoral Lie superalgebras, and in particular compute the tensor product of arbitrary irreducible representations.  Section 4 is a reminder on highest weight irreducible modules over $\fg$.  Section 5 introduces the Grothendieck and reduced Grothendieck rings, and studies the structure of these rings for quasitoral Lie superalgebras; so-called spoiled superalgebras are introduced.  Section 6 presents a general result on the form of the supercharacter of an irreducible highest weight $\fg$-module in the case when the block is nice enough.  Section 7 discusses what more we can say when we assume that $\fg$ is quasireductive and we work with finite-dimensional modules semisimple over $\fg_{\ol{0}}$.  Section 8 introduces the map $ds_x$ on reduced Grothendieck rings, and realizes it as coming from a restriction functor in many cases.  In Section 9 the case of $\fg=\fq_n$ is studied, and we relate the reduced Grothendieck ring of the integral subcategory to a natural quotient of an exterior algebra on an infinite-dimensional vector space. 

\subsection{Acknowledgements} The authors are grateful to I. Entova-Aizenbud, D. Grantcharov, and V. Hinich for numerous helpful discussions.  The first author was supported in part by ISF Grant 1957/21 and the Minerva foundation with funding from the Federal German Ministry for Education and Research.   The second author was supported in part by NSF grant 2001191.  The third author was supported in part by ISF Grant 711/18, NSF-BSF Grant 2019694, and ARC
grant DP210100251

	\subsection{Table of contents} 
\tableofcontents

\subsection{List of notation}

\begin{itemize}
	\item $\cM(\fh)$: see Sec. \ref{section-h-modules}.
	\item $C(\lambda)$: see Sec. \ref{section-irreducible-h-mod}.
	\item $K_{\lambda}$: see Sec. \ref{section-Uh-cliff-alg}.
	\item $F_{\lambda}$: see Sec. \ref{section-Uh-cliff-alg}.
	\item $L(\lambda)$: see Sec. \ref{section-highest-wt-mod}.
	\item $\Gr(\CC)$: see Sec. \ref{section-groth-groups}.
	\item $\Gr_+(\CC)$: see Sec. \ref{section-groth-groups}.
	\item $\Gr_-(\CC)$: see Sec. \ref{section-groth-groups}.
	\item $M_{gr}$, $M_{gr,\pm}$: see Sec. \ref{section-groth-groups}.
	\item $I_0$, $I_1$: see Sec. \ref{section-irreducible-h-mod}.
	\item $\ft_{G}$: see Sec. \ref{section-equivariant-setting}.
	\item $\operatorname{sch}_{\fh}$: see Sec. \ref{section-superchar-map}.
	\item $\Xi$: see Sec. \ref{section-monoid-Xi-cores}.
	\item $\operatorname{Core}$: see Sec. \ref{section-monoid-Xi-cores}.
	\item $\mathcal{F}(\fg)$: see Sec. \ref{section-quasired}.
	\item $P^+(\fg)$: see Sec. 7.
	\item $P(\fg)$: see Sec. 7.
	\item $P(\fg)'$: see Sec. \ref{section-quasired}.
	\item $\operatorname{smult}$: see Sec. \ref{section-ds-rk-s-q(n)}.
	\item $ds_s$: see Sec. \ref{section-ds-rk-s-q(n)}.
\end{itemize}

\section{Preliminaries: maximal quasitoral subalgebras and Clifford algebras}\label{section-prelim}

We work over the field over of complex numbers, $\mathbb{C}$,
and denote by $\mathbb{N}$ the set of non-negative integers.
For a super vector space $V$ we write $V=V_{\ol{0}}\oplus V_{\ol{1}}$ for its parity decomposition.  Then $\Pi V$ will denote the parity-shifted super vector space obtained from $V$.  Let $\delta_V:=\delta$ denote the endomorphism of $V$ given by $\delta(v)=(-1)^{p(v)}v$.

We work with Lie superalgebras which are not necessarily finite-dimensional.  For a Lie superalgebra $\fg$ we denote by $\Fin(\fg)$ the category of finite-dimensional $\fg$-modules.

\subsection{Maximal (quasi)toral subalgebras}\label{section-maxl-qtoral-subalg}

\subsubsection{Definition} Let $\fg$ be Lie superalgebra.  We say that a finite-dimensional subalgebra $\ft\sub\fg_{\ol{0}}$ is a {\em maximal toral} subalgebra if it is commutative, acts diagonally on $\fg$ under the adjoint representation, and we have $\fg_{\ol{0}}^{\ft}=\ft$.  In this case we set $\fh:=\fg^\ft$, and we refer to $\fh$ as a {\em maximal quasitoral} subalgebra of $\fg$.  Observe that $\fh_{\ol{0}}=\ft$.

Denote by $\Delta(\fg):=\Delta\sub\ft^*$ the nonzero eigenvalues of $\ft$ in $\operatorname{Ad}\fg$, and write $Q=\mathbb{Z}\Delta$.  We will assume throughout that
\[
	\text{ all eigenspaces $\fg_{\nu}$ ($\nu\in \Delta\cup\{0\}$) are finite-dimensional.} \ \ \ \ \ \ (*)
\]
In particular we assume that $\fh$ is finite-dimensional.  We have
\[
\fg=\fh\oplus (\bigoplus_{\alpha\in\Delta}\fg_{\alpha}).
\]

\subsubsection{Triangular decomposition}\label{section-triang-decomp}
We choose a group homomorphism $\gamma:\mathbb{Z}\Delta\to\mathbb{R}$ such that $\gamma(\alpha)\neq0$ for all $\alpha\in\Delta$.  Such a homomorphism exists because $\fh\cong\mathbb{R}$ as a $\mathbb{Q}$-vector space.  We introduce the triangular decomposition 
$\Delta(\fg)=\Delta^+(\fg)\coprod \Delta^-(\fg)$,
with
$$\Delta^{\pm}(\fg):=\{\alpha\in\Delta(\fg)|\ \pm\gamma(\alpha)>0\},$$
 and  define a  partial order on $\ft^*$  by
$$\lambda\geq \nu \ \text{ if }\ \nu-\lambda\in\mathbb{N}\Delta^-.$$
We set 
$\fn^{\pm}:=\oplus_{\alpha\in\Delta^{\pm}}\fg_{\alpha}$ and call a subalgebra of the form $\fb:=\fh\oplus\fn^+$ a Borel subalgebra.  We further assume throughout that
\[
\mathcal{U}(\fn^-)_{\nu} \text{ is finite-dimensional for all }\nu\in Q\ \ \ \ \ \ (**)
\]

\begin{remark}
	Several notions of triangular decompositions for Lie superalgebras have appeared in the literature. In \cite{M} and \cite{PS} for example, a notion of positive roots arose from the choice of a generic hyperplane in $\fh^*$.  Our approach generalizes these approaches and admits a more flexible definition. One can construct finite-dimensional Lie superalgebras for which our definition gives rise to more triangular decompositions as compared to \cite{M} and \cite{PS}, e.g.~consider $\fg$ with $\fg_{\ol{0}}$ one-dimensional acting by real, $\mathbb{Q}$-linearly independent characters on $\fg_{\ol{1}}$.
	
	On the other hand, for simple, finite-dimensional Lie superalgebras	our notion of triangular decomposition agrees with that of \cite{PS}.
\end{remark}

\subsubsection{Examples}
\begin{itemize}
	\item If $\fg$ is a Kac-Moody superalgebra, then, by~\cite{KP} any maximal toral subalgebra $\ft$ satisfying
$(*)$ is Cartan subalgebra of $\fg_{\ol{0}}$; one has $\fh=\ft$.
	\item If  $\fg$ is a quasi-reductive Lie superalgebra 
	($\dim \fg<\infty$, $\fg_{\ol{0}}$ is reductive and
	$\fg_{\ol{1}}$ is a semisimple $\fg$-module), then a maximal toral (resp. quasitoral)  subalgebra
	$\ft$ is a Cartan subalgebra of $\fg_{\ol{0}}$ (resp. $\fg$).    In both this example and the former, $\ft$ and $\fh$ are unique
	up to a conjugation by inner automorphism, see~\cite{KP}.

	\item If we fix an invariant form on a quasi-reductive Lie superalgebra $\fg$
	(this can be the zero form), we can construct the affinization
	$\hat{\fg}$ with  $\hat{\ft}=\ft+\mathbb{C}K+\mathbb{C}d$
	and $\hat{\fh}=\fh+\mathbb{C}K+\mathbb{C}d$.
		
	\item The cases when $\ft\not=\fh$ include the queer Lie superalgebras and their affinizations.
\end{itemize}

\subsection{Clifford algebras}

For a vector space $V$ with a symmetric (not necessarily nondegenerate) bilinear form $F$, let $\Cl(V,F)$ denote the corresponding Clifford algebra.  This is a superalgebra where the elements of $V$ are declared to be odd.  Write $K\sub V$ for the kernel of $F$, so that $F$ induces a nondegenerate form on $V/K$, which we also write as $F$.  Then we have an isomorphism of superalgebras:
\[
\Cl(V,F)\cong\Cl(V/K,F)\otimes\Lambda^\bullet K.
\]
The superalgebra $\Cl(V,F)$ is semisimple if and only if $F$ is nondegenerate.  

\subsubsection{}
Suppose $F$ is nondegenerate, and write $m$ for the dimension of $V$.  We have
\[
\Cl(V,F)\cong Q(2^{n+1})\ \ \text{ if }m=2n+1, \ \ \Cl(V,F)\cong 
\End(\mathbb{C}^{2^{n-1}|2^{n-1}} )\ \ \text{ if }m=2n.
\]
Here $Q(r)$ is the queer superalgebra, i.e. the associative subalgebra of $\End(\mathbb{C}^{r}|\mathbb{C}^{r})$
consisting of the matrices of the form $\begin{pmatrix} A& B\\ B &A\end{pmatrix}$.

  It follows that $\Cl(V,F)$ is a simple superalgebra and admits a unique, parity invariant, irreducible module when $m$ is odd, while if $m$ is even there are two irreducible modules that differ by parity.   Moreover, all ($\mathbb{Z}_2$-graded) $\Cl(V,F)$-modules are completely reducible. For $m\not=0$, if $E$ is an irreducible representation of $\Cl(V,F)$ one has 
\[
\dim E_{\ol{0}}=\dim E_{\ol{1}}=2^{\lfloor\frac{m-1}{2}\rfloor}.
\]

\subsubsection{} 
For a non-negative integer $m$, we write $\Cl(m)$ for the Clifford algebra $\Cl(\mathbb{C}^m,F)$, where $F$ is the standard nondegenerate symmetric bilinear form on $\mathbb{C}^m$. Clearly $\Cl(V,F)\cong\Cl(m)$ if $\dim V=m$ and $F$ is nondegenerate.

\subsubsection{}\label{section-outer-tensor}
Let $A_1,A_2$ be associative superalgebras and $V_i$ be $A_i$-modules;
we define the outer tensor product  $V_1\boxtimes V_2$ as 
the space $V_1\otimes V_2$ endowed by the $A_1\otimes A_2$-action
$$(a_1,a_2)(v_1\otimes v_2):=(-1)^{p(a_2)p(v_1)} a_1v_1\otimes a_2v_2.$$

One has $\Cl(m)\otimes \Cl(n)\cong\Cl(m+n)$; 
if  $V_1$ (resp., $V_2$) are simple modules over $\Cl(m_1)$ (resp., 
$\Cl(m_2)$), then $V_1\boxtimes V_2$ is a simple if either $m_1$ or $m_2$ is even, and if $m_1$ and $m_2$ are both odd then $V_1\otimes V_2=L\oplus\Pi L$, where $L$ is simple over $\Cl(n_1+m_1)$.

\subsection{Realization of the irreducible representation of $\Cl(2n)$}\label{section-clif-realization}
Consider $\mathbb{C}^{2n}$ with standard basis $e_1,\dots,e_n$, $f_1,\dots,f_n$, equipped with the symmetric nondegenerate form $(-,-)$ satisfying:
\[
(e_i,f_j)=\delta_{ij},\ \ \ \ \ \ (e_i,e_j)=(f_i,f_j)=0.
\] 
Consider the polynomial superalgebra $L=\mathbb{C}[\xi_1,\dots,\xi_n]$ with odd generators $\xi_1,\dots,\xi_n$. Then we may realize $L$ as an irreducible representation of $\Cl(2n)$ via $e_i\mapsto\xi_i$ and $f_j\mapsto\d_{\xi_j}$, i.e. $e_i$ acts by multiplication by $\xi_i$ and $f_j$ acts by the derivation sending $\xi_{i}\mapsto\delta_{ij}$.  In this way we have defined a surjective morphism $\Cl(2n)\to \End(L)$ (in fact it is an isomorphism).  Every irreducible representation of $\Cl(2n)$ is isomorphic to $L$ or $\Pi L$.

\subsubsection{} Continuing with the setup from  Section \ref{section-clif-realization}, if we choose $W$ a maximal isotropic subspace of $\mathbb{C}^{2n}$, then $\Lambda^\bullet W$ acts on $L$, and under this action $L$ is isomorphic to the exterior algebra of $W$ under left multiplication. Thus given two arbitrary irreducible representations of $\Cl(V)$, they are isomorphic if and only if the parities of the one-dimensional $W$-invariant subspaces of each are the same. 

\subsection{The case $\Cl(2n+1)$}
Next we look at $\mathbb{C}^{2n+1}$ with nondegenerate symmetric form $(-,-)$ and basis $e_1,\dots,e_n$, $f_1,\dots,f_n$, $g$.  Here the inner product relations for the $e_i$'s and $f_i$'s are the same as in Section \ref{section-clif-realization}, along with
\[
(g,e_i)=(g,f_j)=0, \ \ \ \ (g,g)=2.
\]
Consider the exterior algebra $L=\mathbb{C}[\xi_1,\dots,\xi_{n+1}]$.  Then we may realize $L$ as the unique irreducible representation of $\Cl(2n+1)$ by $e_i\mapsto\xi_i$, $f_j\mapsto\d_{\xi_j}$, and $g\mapsto \xi_{n+1}+\d_{\xi_{n+1}}$. 

\subsection{The operator $T$}\label{section-cliff-T}

Choose an orthonormal basis $H_1,\dots,H_m$ of $\mathbb{C}^m$, and let \newline $T=H_1\cdots H_m\in\Cl(m)$. Note that this operator is an eigenvector of $O(m)$ with weight given by the determinant.  Thus it is well-defined up to an orientation on $V$.  

Let $V$ be an irreducible representation of $\Cl(m)$.  If $m$ is odd, then $\Pi V\cong V$, and $\End_{\mathbb{C}}(V)=\Cl(m)\oplus \Cl(m)\delta_V$ (here we consider all endomorphisms), where $\delta_V(v)=(-1)^{p(v)}v$. In this case, $T=\phi\delta_V$, where $\phi$ is an odd $\Cl(m)$-equivariant automorphism of $V$.  If $m$ is even, then $V\not\cong\Pi(V)$ and $\Cl(m)=\End(V)$.  In this case, $T=(-1)^{n}\delta_V\in\End(V)$, where $\dim V=2n$.

\section{Representation theory of quasitoral Lie superalgebras} 

\subsection{$\fh$-modules}\label{section-h-modules}
Take $\fh$ as in Section \ref{section-maxl-qtoral-subalg}:
$\fh$ is a finite-dimensional Lie superalgebra with 
$$[\ft,\fh]=0,$$
where $\ft=\fh_{\ol{0}}$.  We call Lie superalgebras of this form quasitoral.  For  a semisimple $\ft$-module $N$, and $\nu\in\ft^*$, write $N_{\nu}$ for the $\nu$-weight space in $N$.  

Denote by $\cM(\fh)$ the full subcategory of $\fh$-modules $N$ with diagonal action of $\ft$ and finite-dimensional weight spaces
$N_{\nu}$.  We set $\cF(\fh)$ to be the full subcategory of $\cM(\fh)$ consisting of those modules which are finite-dimensional.  The simple modules in $\cM(\fh)$ and $\cF(\fh)$ coincide.  In this section we study the questions of restriction, tensor product, and extensions of simple modules in $\cM(\fh)$. 

We denote by $\sigma$ the anti-automorphism of $\cU(\fh)$ induced by the antipode
$-\Id|_{\fh}$ (recall that anti-automorphism means 
$\sigma(ab)=(-1)^{p(a)p(b)}\sigma(b)\sigma(a)$).
 This anti-automorphism induces the standard duality $*$ on $\cF(\fh)$.

\subsubsection{}\label{section-Uh-cliff-alg}  View  $\cU(\fh)$  as a Clifford algebra over the polynomial algebra $\cS(\ft)$; the corresponding symmetric bilinear form is
given on $\fh_{\ol{1}}$ by the formula $F(H,H')=[H,H']$
(see Appendix in~\cite{Gq} for details). 

For each $\lambda\in\ft^*$, the evaluation of $F$ at $\lambda$ gives a symmetric form 
$F_{\lambda}: (H,H')\mapsto\lambda([H,H'])$.  We denote by $\operatorname{rk}F_{\lambda}$ the rank of this form.
 For each $\lambda\in\ft^*$ we consider the Clifford algebra
\[
\Cl(\lambda):=\Cl(\fh_{\ol{1}},F_{\lambda})=\cU(\fh)/\cU(\fh) I(\lambda),
\]
 where
$I(\lambda)$ stands for the kernel
of the algebra homomorphism $\cS(\ft)\to\mathbb{C}$ induced by $\lambda$.  We will write $K_{\lambda}\sub\fh_{\ol{1}}$ for the kernel of $F_{\lambda}$.  
Then we have an isomorphism of superalgebras
\begin{equation}\label{eqn-eqcliff}
\Cl(\lambda)\cong \Cl(\rank F_{\lambda})\otimes \bigwedge K_{\lambda}.
\end{equation}

Denote by $\phi_{\lambda}: \cU(\fh)\to \Cl(\lambda)$  the canonical
epimorphism, and $p_{\lambda}:\Cl(\lambda)\to\Cl(\lambda)/(K_{\lambda})$ for the projection. Let $\varphi:\fh\to\fh$ be an automorphism of Lie superalgebras, and write also $\varphi$ for the corresponding automorphism of $\cU(\fh)$.  Then for every $\lambda$, $\varphi$ induces isomorphisms of algebras
\[
\varphi_{\lambda}:\Cl(\lambda)\to\Cl(\varphi^{-1}(\lambda)), \ \ \ \ \ol{\varphi_{\lambda}}:\Cl(\lambda)/K_{\lambda}\to\Cl(\varphi^{-1}(\lambda))/(K_{\varphi^{-1}(\lambda)}).
\]
We have the following commutative diagram:
\begin{equation}
\xymatrix{\cU(\fh) \ar[d]^{\phi_{\lambda}} \ar[rr]^{\varphi} & & \cU(\fh)\ar[d]^{\phi_{\varphi^{-1}(\lambda)}}  \\
           \Cl(\lambda)  \ar[rr]^{\varphi_{\lambda}} \ar[d]^{p_{\lambda}} &  &  \Cl(\varphi^{-1}(\lambda)) \ar[d]^{p_{\varphi^{-1}(\lambda)}} \\
       \Cl(\lambda)/(K_{\lambda}) \ar[rr]^{\ol{\varphi_{\lambda}}} & & \Cl(\varphi^{-1}(\lambda))/(K_{\varphi^{-1}(\lambda)}) }
\end{equation}
For the anti-involution $\sigma$ we also have the same diagram as above, where the induced maps $\sigma_{\lambda},\ol{\sigma_{\lambda}}$  are anti-algebra isomorphisms.

\begin{lemma}\label{lemma-cl(lambda)-modules}
\begin{enumerate}
	\item A $\Cl(\lambda)$-module $N$ is semisimple if and only if
	$\Ann_{\fh_{\ol{1}}} N=K_{\lambda}$.
	
	\item If $N$ is an indecomposable $\Cl(\lambda)$-module of length $2$, 
	then $[N:C(\lambda)]=[N:\Pi(C(\lambda))]$.
\end{enumerate}
\end{lemma}

\begin{proof}
	These follow from Formula~(\ref{eqn-eqcliff}).
\end{proof}

% \subsubsection{} We say the form $F$ on $\fh_{\ol{1}}$ is diagonalizable if there exists a basis $H_1,\dots,H_n$ for $\fh_{\ol{1}}$ such that $F(H_i,H_j)=0$ for $i\neq j$.

\subsubsection{Examples}\begin{itemize}
	\item If $\fh$ is quasitoral such that $\fh_{\ol{1}}$ is commutative, then $F$ is the zero form.
	
	\item For $\fg=\fq_n$ and $\fh$ a maximal quasitoral subalgebra, one has $\dim\fh=(n|n)$ and $\fh\cong\fq_1\times\cdots\times\fq_1$.  Thus in this case the form $F$ is diagonal, and $\operatorname{rk}F_{\lambda}$ is the number of nonzero entries of $\lambda$ under the decomposition $\ft\cong(\fq_1)_{\ol{0}}\times\cdots\times(\fq_1)_{\ol{0}}$.

	\item For $\fg=\fsq_n$ and $\fh$ a maximal quasitoral subalgebra, one has $\dim\fh=(n|n-1)$.  In this case $F$ is not diagonal.
\end{itemize}

\subsection{Irreducible $\fh$-modules}\label{section-irreducible-h-mod}
The irreducible $\fh$-modules all arise from irreducible modules over $\Cl(\lambda)$ for some $\lambda\in\ft^*$.  We denote by $C(\lambda)$ a simple $\Cl(\lambda)$-module and also view it as an $\fh$-module. For $\lambda=0$ we fix the grading by taking $C_0=\mathbb{C}$; for all
other values of $\lambda$ we fix a grading in an arbitrary way until further notice.  By above, 
\begin{equation}\label{eqn-C-lambda}\begin{array}{l}
\dim C(\lambda)=2^{n_{\lambda}}\ \text{ where }
n_{\lambda}:=\lfloor\frac{\rk F_{\lambda}+1}{2}\rfloor,\\
\{u\in \fh_{\ol{1}}|\ uC(\lambda)=0\}=K_{\lambda}.\end{array}
\end{equation}

Set:
\begin{equation}\label{eqn-I01}
I_i=\{\lambda\in\ft^*: \rk F_{\lambda}\equiv i \mod 2\} \ \text{ for }i=0,1.
\end{equation}
 Then $C(\lambda)\cong \Pi C(\lambda)$ if and only if  $\lambda\in I_1$;
we will often use the notation $\Pi^{\frac{\rank F_{\lambda}}{2}} C(\lambda)$, where $\Pi^{\frac{\rank F_{\lambda}}{2}} C(\lambda)\cong C(\lambda)\cong 
\Pi C(\lambda)$ whenever $\rank F_{\lambda}$ is odd, and if $\rank F_{\lambda}$ is even it has the obvious meaning.

\subsection{Blocks of $\cM(\fh)$ and $\cF(\fh)$} The blocks of both $\cM(\fh)$ and $\cF(\fh)$ are parameterized up to parity shift by $\lambda\in\ft^*$, as follows:  If $\corank F_{\lambda}>0$, then there is one block of both $\cM(\fh)$ and $\cF(\fh)$ on which $\ft$ acts by the character $\lambda$; in both cases this block is equivalent to the category of finite-dimensional modules over $\bigwedge K_{\lambda}$.  If $\corank F_{\lambda}=0$ then for both categories there is one (resp. two) block(s) of $\fh$ on which $\ft$ acts by $\lambda$ when $\lambda\in I_1$ (resp. $\lambda\in I_0$), and the block(s) is (are) semisimple.  

\subsubsection{Remark}  One can think of $\corank F_{\lambda}$ as the ``atypicality'' of its corresponding block.  In particular $\corank F_{\lambda}=0$ if and only if the block is semisimple and thus its objects are projective in $\cM(\fh)$, and in general the block corresponding to $\lambda$ is equivalent to modules over a Grassmann algebra on $\text{corank}(F_{\lambda})$-many variables.

\subsection{The operator $T_{\fh}$} Let $H_1,\dots,H_n$ be a basis of $\fh_{\ol{1}}$, and define 
\[
T_{\fh}=\{H_1,\{\cdots \{H_n,1\}\cdots\}\in\cU(\fh)
\]  
where $\{x,y\}=xy+(-1)^{\ol{x}\ol{y}}yx$ denotes the super anticommutator in $\cU(\fh)$.  It is known that up to scalar, $T_{\fh}$ does not depend on the choice of basis, see~\cite{Gg}. This operator anticommutes with $\fh_{\ol{1}}$, so in particular the image of a submodule under $T_{\fh}$ remains a submodule.

\subsubsection{Action of $T_{\fh}$ on simples} The action of $T_{\fh}$ on simple $\fh$-modules is deduced from Section \ref{section-cliff-T}, and is as follows.  If $\corank F_{\lambda}>0$, then $T_{\fh}$ acts by $0$.  If $\corank F_{\lambda}=0$, then $T_{\fh}$ acts by an automorphism, although not $\fh$-equivariantly.  If $n$ is odd, then $T_{\fh}$ is a nonzero multiple of $\delta\phi$, where $\phi$ is an $\fh$-equivariant odd automorphism.  If $n$ is even, $T_{\fh}$ is a nonzero multiple of $\delta$.

In particular when $n$ is even, $T$ acts on $C(\lambda)$ by an operator of the form 
\[
a(\lambda)\Id_{(C(\lambda))_{\ol{0}}}\oplus(-a(\lambda))\Id_{(C(\lambda))_{\ol{1}}},
\]
for a scalar $a(\lambda)$. Thus $T$ distinguishes between $C(\lambda)$ and its parity shift for projective irreducible modules.  

\subsubsection{Action of $T_{\fh}$ on all of $\cM(\fh)$}  Let $\corank F_{\lambda}>0$.  Then the injective hull of $C(\lambda)$ is given by the $\Cl(\lambda)$-module $I(C(\lambda))=C(\lambda)\otimes\Lambda^\bullet K_{\lambda}$.  We claim that 
\begin{enumerate}
	\item $T_{\fh}$ annihilates the radical of $I(C(\lambda))$;
	\item $\operatorname{Im}T_{\fh}=C(\lambda)=\operatorname{socle}I(C(\lambda))$.
\end{enumerate}
In other words, $T_{\fh}$ acts by taking the head of this module to its socle.  It follows that we understand completely the action of $T_{\fh}$ on every module in $\cM(\fh)$.  

To prove our claim, choose a basis $f_1,\dots,f_r$ of $K_{\lambda}$ and extend it to a basis $e_1,\dots,e_s,f_1,\dots,f_r$ of $\fh_{\ol{1}}$ so that $F_{\lambda}(e_i,e_{j})=\delta_{ij}$.  Then the image of $T_{\fh}$ in $\Cl(\fh)$, up to scalar, is given by
\[
e_1\cdots e_sf_1\cdots f_r.
\]
By considering the action of this operator, the statement is clear.

\subsubsection{Remark}\label{remark-rigidify-parity} 
It is possible to uniformly choose the parity of irreducible $\fh$-modules as follows.  
We consider a linear order on $\mathbb{C}$
given by $c_2>c_1$ if the real part of $c_2-c_1$ is positive or the real part is zero and the imaginary part is positive.
We fix any function 
$$t: \ft^*\to\{c\in\mathbb{C}|\ c>0\}.$$
Retain notation of Section \ref{section-Uh-cliff-alg}.
For each $\lambda$ 
we choose $T_{\lambda}\in \cU(\fh)$ in such a way that 
\[
\ol{T}_{\lambda}=p_{\lambda}\circ\phi_\lambda(T_{\lambda})\in \Cl(\lambda)/K_{\lambda}
\] 
is an anticentral element satisfying $\ol{T}_{\lambda}^2=t(\lambda)^2$
(the element $\ol{T}_{\lambda}$ is unique up to sign).

If $\lambda\in I_0$, then $T_{\lambda}$ is even 
and it  acts
on $C(\lambda)$ by a nonzero superconstant $\pm t(\lambda)$ and
we fix a grading on $C(\lambda)$ 
by taking 
$$(C(\lambda))_{\ol{0}}:=\{v\in C(\lambda)|\ T_{\lambda}v=t(\lambda)v\}.$$

\subsubsection{Dualities in $\cF(\fh)$}\label{section-dualities-fd}
Fix $\lambda\in \ft^*$. The category $\cF(\fh)$ has the duality $*$ induced by 
$\sigma$, and
another contragredient involution
$(-)^\#: \cF(\fh)\to\cF(\fh)$  induced by the antiautomorphism $\sigma'(a)=a$ for $a\in\fh_{\ol{0}}$ and $\sigma'(a)=\sqrt{-1}a$ for $a\in\fh_{\ol{1}}$.
Note that $\sigma'$ induces an anti-involution on $\Cl(\lambda)$.

The element $\ol{T}_{\lambda}$ can be written as the product $H_1'\ldots H'_k$, where
$H'_1,\ldots,H_k'$ is a lift of a basis
of  $\Cl(\lambda)/K_{\lambda}$ 
satisfying $[H'_i,H'_j]=0$ for $i\not=j$. Therefore
 $\ol{T}_{\lambda}=(-1)^{\frac{\rank F_{\lambda}}{2}} \sigma'(\ol{T}_{\lambda})$ 
for $\lambda\in I_0$; this gives the following useful formula
\begin{equation}\label{eqn-frisk}
C(\lambda)^{\#}\cong \Pi^{\frac{\rk F_{\lambda}}{2}}
C(\lambda);\end{equation}
which was first established in~\cite{Frisk}, Lemma 7.   

Since
$\ol{\sigma_{\lambda}}:\Cl(\lambda)/K_{\lambda}\to \Cl(-\lambda)/K_{-\lambda}$ is anti-isomorphism
we have
$\frac{\ol{\sigma_{\lambda}}(\ol{T}_{\lambda})}{t(\lambda)}=(-1)^i 
\frac{\ol{T}_{-\lambda}}{t(-\lambda)}$ for some 
$i\in \{0,1\}$, and correspondingly $C(\lambda)^{*}\cong \Pi^{i} C(-\lambda)$.

\subsection{Restriction to quasitoral subalgebra}  Given a quasitoral subalgebra $\fh'\sub\fh$, we have $\ft'=\fh'_{\ol{0}}\sub\fh_{\ol{0}}=\ft$.  Thus we have a natural restriction $\ft^*\to(\ft')^*$, and therefore we consider weights $\lambda\in\ft^*$ as defining weights in $(\ft')^*$ naturally.  We write $F_{\lambda}'$ for the bilinear form induced on $\fh'_{\ol{1}}$ by a given weight $\lambda\in\ft^*$, which is exactly the restriction of $F_{\lambda}$ to $\fh_{\ol{1}}'$.

Let $\Cl'(\lambda)$ be the subalgebra of $\Cl(\lambda)$ which is generated by $\fh_{\ol{1}}'$;
clearly, $\Cl'(\lambda)=\Cl(\fh_{\ol{1}}',F'_{\lambda})$.  Denote by $E'({\lambda})$ a simple $\Cl'(\lambda)$-module.

\begin{proposition}\label{proposition-cliff-res}
Write $V'=\fh_{\ol{1}}'$ and $V=\fh_{\ol{1}}$.
\begin{enumerate}
\item $C(\lambda)$ is simple over $\Cl'(\lambda)$ if and only if 
$\lfloor\frac{\rk F_{\lambda}+1}{2}\rfloor=\lfloor\frac{\rk F'_{\lambda}+1}{2}\rfloor$.

\item $C(\lambda)$ is 
semisimple over $\Cl'(\lambda)$ if and only if 
$\Ker F'_{\lambda}=V'\cap \Ker F_{\lambda}$.

\item If $\rk F_{\lambda}\not=\rk F'_{\lambda}$, then 
$[C(\lambda):E'_{\lambda}]=[C(\lambda):\Pi E'_{\lambda}]$.
\end{enumerate}
\end{proposition}

\begin{proof}
The case (i) follows 
from~(\ref{eqn-C-lambda}) and case (ii)
 follows from part (ii) of lemma \ref{lemma-cl(lambda)-modules}.

 For (iii)  assume that  $\rk F_{\lambda}\not=\rk F'_{\lambda}$.
Substituting $\fh$ by $\fh/\Ker F_{\lambda}$
we may assume that $\Ker F_{\lambda}=0\not=\Ker F'_{\lambda}$.
Then $\fh_{\ol{1}}$ admits a basis $H_1',\ldots, H'_m$ such that
$V'$ is spanned by $H_1,\ldots, H'_p$ and
the matrix of $F_{\lambda}$ takes the form
$$\begin{pmatrix}
Id_{p-s} & 0 & 0& 0\\
0& 0& Id_s & 0\\ 
0& Id_s& 0 & 0 \\
0 & 0& 0 & Id_k
\end{pmatrix}$$
with $k+s+p=m$, $k,s\geq 0$ and $k+s>0$ (since $V'\not=V$).
Note that $\Ker F'_{\lambda}$ is spanned by $H_{p+1},\ldots, H'_{p+s}$
and so $E'_{\lambda}$ is a simple $\Cl(p-s)$-module.
If $k\not=0$, then 
 the action of $H'_k$ to $C(\lambda)$
is an odd involutive $\Cl'(\lambda)$-homomorphism, so $\Res_{\Cl'(\lambda)} C(\lambda)$
is $\Pi$-invariant. Consider the remaining case $k=0$.   
Then $s\not=0$ and  $\Cl(\lambda)=\Cl(p-s)\otimes \Cl(2s)$.
Using~\ref{section-outer-tensor} we get
$C(\lambda)\cong E'_{\lambda}\boxtimes E''$, 
where $E''$ is a simple $\Cl(2s)$-module. By above, 
$\dim E''_{\ol{0}}=\dim E''_{\ol{1}}$. Hence 
$\Res_{\Cl(p-s)} C(\lambda)$
is $\Pi$-invariant, that is
$$[C(\lambda): E'_{\lambda}]=[C(\lambda): \Pi E'_{\lambda}]$$
as required. This establishes (iii).
\end{proof}

\subsection{Tensor product of irreducible $\fh$-modules}

Let $\lambda,\mu\in\ft^*$.  We compute \newline $C(\lambda)\otimes C(\mu)$.  Observe that $C(\lambda)\otimes C(\mu)$ is naturally a module over $\Cl(\fh_{\ol{1}}/K_{\lambda}\cap K_{\mu},F_{\lambda+\mu})$ and $K_{\lambda}\cap K_{\mu}\subseteq K_{\lambda+\mu}$.  Set
\[
K_{\lambda,\mu}:=K_{\lambda+\mu}/(K_{\lambda}\cap K_{\mu}).
\]

We have an isomorphism of superalgebras
\[
\Cl(\fh_{\ol{1}}/K_{\lambda}\cap K_{\mu},F_{\lambda+\mu})\cong \Cl(\fh_{\ol{1}}/K_{\lambda+\mu},F_{\lambda+\mu})\otimes\bigwedge K_{\lambda,\mu}.
\]

\begin{lemma}
$C(\lambda)\otimes C(\mu)$ is projective over $\Cl(\fh_{\ol{1}}/(K_{\lambda}\cap K_{\mu}),F_{\lambda+\mu})$.
\end{lemma}

\begin{proof}
	It suffices to show that $\bigwedge K_{\lambda,\mu}$ acts freely.  Let $v\in K_{\lambda,\mu}$.  Then WLOG $v\notin K_{\lambda}$, so the subalgebra generated by $v$ acts projectively on $C(\lambda)$, and thus also on the tensor product.  The statement now follows from facts about the representation theory of exterior algebras (see \cite{AAH}). 
\end{proof}

\subsubsection{}
Notice that the unique (up to parity) indecomposable projective module $P$ over $\bigwedge K_{\lambda,\mu}$ is the free module of rank 1.  Thus we have shown that $C(\lambda)\otimes C(\mu)$ is a sum of modules of the form $(\Pi)C(\lambda+\mu)\otimes\bigwedge K_{\lambda,\mu}$.  

If the rank of either $F_{\lambda}$ or $F_{\mu}$ is odd, or the rank of $F_{\lambda+\mu}$ is odd, then the tensor product $C(\lambda+\mu)\otimes P$ is parity invariant, so the explicit decomposition of $C(\lambda)\otimes C(\mu)$ is the appropriate number of copies of $C(\lambda+\mu)\otimes P$ and its parity shift, according to a dimension count.

\subsubsection{}
Thus let us suppose that $\rk F_{\lambda},\rk F_{\mu},$ and $\rk F_{\lambda+\mu}$ are all even and we have $\rk F_{\nu}=2n_{\nu}$ for $\nu=\lambda,\mu,\lambda+\mu$.  

By Section \ref{section-clif-realization}, we may realize $C(\lambda)$ as $k[\xi_1,\dots,\xi_n]$ and $C(\mu)$ as $k[\eta_1,\dots,\eta_m]$, so that 
\[
C(\lambda)\otimes C(\mu)=k[\xi_1,\dots,\xi_n,\eta_1,\dots,\eta_m].
\]
Now choose a maximal isotropic subspace $U$ for $F_{\lambda+\mu}$ in $\fh_{\ol{1}}$.  Then the parity decomposition of $C(\lambda)\otimes C(\mu)$ is described by the $U$-invariants on this space.  Let $u\in U$.  Then because $u^2$ acts trivially on this module, its action on the tensor product is given by (using Section \ref{section-clif-realization})
\[
u\mapsto \sum\limits_{i=1}^na_iX_i+\sum\limits_{j=1}^mb_jY_j,
\]
where $a_i,b_j\in\mathbb{C}$, and $X_i\in\{\xi_i,\d_{\xi_i}\}$ and $Y_j\in\{\eta_j,\d_{\eta_j}\}$.  Further, because $[U,U]$ acts trivially, we can choose $X_1,\dots,X_n$ and $Y_1,\dots,Y_m$ uniformly so that every element of $U$ acts in the way described for $u$, with potentially different coefficients $a_i$ and $b_j$.  Now, suppose that $X_i=\xi_i$ for some $i$.  Define an odd, linear automorphism $s_i$ of $k[\xi_1,\dots,\xi_n]$ as follows.  For $J=\{i_1,\dots,i_{|J|}\}\sub\{1,\dots,n\}\setminus\{i\}$, write $\xi_{J}=\xi_{i_1}\cdots\xi_{i_{|j|}}$, and then set:
\[
s_i(\xi_J)=\xi_i\xi_J \ \ \ \ s_i(\xi_i\xi_J)=\xi_J.
\]
Then under this automorphism, multiplication by $\xi_i$ becomes $\d_{\xi_i}$ and vice versa, while for $i\neq j$, $\d_{\xi_j}$ and multiplication by $\xi_j$ become negative themselves.  Using this automorphism, we may instead assume that $X_i=\d_{\xi_i}$, and in this way we may assume that $X_i=\d_{\xi_i}$ for all $i$, and $Y_j=\d_{\eta_j}$ for all $j$.  Thus $U$ acts by a subspace of constant coefficient vector fields on $k[\xi_1,\dots,\xi_n,\eta_1,\dots,\eta_m]$.  Write $Z\subseteq\langle\xi_1,\dots,\xi_n,\eta_1,\dots,\eta_m\rangle$ for the invariants of $U$ in this subspace.  Then
\[
(C(\lambda)\otimes C(\mu))^{U}={\bigwedge}^\bullet Z.
\]
Thus we have shown the following (recall that $n_{\lambda}:=\lfloor\frac{\rk F_{\lambda}+1}{2}\rfloor$):

\begin{theorem}\label{theorem-tensor-prod}
If $n_{\lambda+\mu}+\dim K_{\lambda,\mu}=n_{\lambda}+n_{\mu}$, then up to parity $C(\lambda)\otimes C(\mu)\cong C(\lambda+\mu)\otimes \bigwedge K_{\lambda,\mu}$.  Otherwise
	\[
	C(\lambda)\otimes C(\mu)=(C(\lambda+\mu)\otimes\bigwedge K_{\lambda,\mu})\otimes \mathbb{C}^{2^a|2^a}.
	\]
	where $a=\frac{n_{\lambda}+n_{\mu}-n_{\lambda+\mu}}{2}-\dim K_{\lambda,\mu}$.
\end{theorem}

\begin{corollary}\label{corollary-tensor-prod-parity-invce}
The module $C(\lambda)\otimes C(\mu)$ is $\Pi$-invariant except for the case when
$K_{\lambda},K_{\mu},K_{\lambda+\mu}$ have even codimensions in $\fh_{\ol{1}}$ and
$$\fh_{\ol{1}}=K_{\lambda}+K_{\mu}+K_{\lambda+\mu}.$$
\end{corollary}

\begin{proof}
Note that $\Pi(C(\lambda))=C(\lambda)\otimes \Pi(\mathbb{C})\cong C(\lambda)$ implies
$\Pi(C(\lambda)\otimes C(\mu))\cong C(\lambda)\otimes C(\mu)$. On the other hand,
 if $\codim\Ker F_{\lambda+\mu}$ is odd, then   any $\fg_{\ol{1}}/ K_{\lambda+\mu}$-module is $\Pi$-invariant. Therefore $C(\lambda)\otimes C(\mu)$ is $\Pi$-invariant if at least one of the numbers $\codim K_{\lambda},\codim K_{\mu},\codim K_{\lambda+\mu}$ is odd. Now assume that these numbers are even. Note  that
$$K_{\lambda}\cap K_{\mu}=K_{\lambda+\mu}\cap K_{\mu}=K_{\lambda}\cap K_{\lambda+\mu}.$$
and set
$$ m_{\lambda,\mu}:=\dim (K_{\lambda}\cap K_{\mu}),\ \ 
r_{\lambda}:=\dim K_{\lambda}/(K_{\lambda}\cap K_{\mu}),\ \ 
r_{\mu}:=\dim K_{\mu}/(K_{\lambda}\cap K_{\mu}).$$

Assume that $C(\lambda)\otimes C(\mu)$ is not $\Pi$-invariant. 
By Theorem \ref{theorem-tensor-prod} in this case
$$n_{\lambda+\mu}+\dim K_{\lambda,\mu}=n_{\lambda}+n_{\mu}.$$
One has
$$n_{\lambda}=\lfloor\frac{\codim K_{\lambda}+1}{2}\rfloor=\frac{\codim K_{\lambda}}{2}=\frac{\dim\fh_{\ol{1}}-r_{\lambda}-m_{\lambda,\mu}}{2}$$
with the similar formulae for $n_{\mu}$ and 
$n_{\lambda+\mu}=\frac{n-\dim K_{\lambda,\mu}-m_{\lambda,\mu}}{2}$.
This gives
$$\dim K_{\lambda,\mu}+m_{\lambda,\mu}+r_{\lambda}+r_{\mu}
=\dim\fh_{\ol{1}}$$
as required.
\end{proof}

\section{The irreducible modules $L(\lambda)$ of $\fg$}

We now return to the setting of Section \ref{section-maxl-qtoral-subalg}, so that $\fg$ denotes a Lie superalgebra containing a finite-dimensional quasitoral subalgebra $\fh$.  Choose a triangular decomposition $\fg=\fn^{-}\oplus\fh\oplus\fn$ as in Section \ref{section-triang-decomp}.

\subsection{Highest weight modules}\label{section-highest-wt-mod}  We call a $\fg$-module $N$ a module of highest weight $\lambda$
if $N_{\lambda}\not=0$ and $N_{\nu}\not=0$ implies $\nu\leq \lambda$.

View $C(\lambda)$ as a $\fb$-module with the zero action of $\fn$
and set
$$M(\lambda):=\Ind^{\fg}_{\fb} C(\lambda);$$  
the module $M(\lambda)$
has a  unique simple quotient which we denote by 
$L(\lambda)$. Each simple module of  highest weight $\lambda$
is isomorphic to $L(\lambda)$ if $\rk F_{\lambda}$ is odd
(i.e., if $\lambda\in I_1$); if $\rk F_{\lambda}$ is even
(i.e., if $\lambda\in I_0$),  each simple module of highest weight $\lambda$
is isomorphic to either  $L(\lambda)$ or to $\Pi L(\lambda)$,
and these modules are not isomorphic.

Note that $\ft$ acts diagonally on $M(\lambda)$ and
all weight spaces $M(\lambda)_{\nu}$ are finite-dimensional
(since we assume all weight spaces $\mathcal{U}(\fn^-)_{\nu}$ are finite-dimensional); in particular, $\Res^{\fg}_{\fh} M(\lambda)$ lies
in $\cM(\fh)$.

\subsection{Duality}\label{section-duality}
In many cases the antiautomorphism $(-)^\#$ introduced in Section \ref{section-dualities-fd} can be extended to an antiautomorphism
of $\fg$ which satisfies $(a^{\#})^{\#}=(-1)^{p(a)}a$.
Using this antiautomorphism we can introduce  a contragredient 
duality on $\fg$-modules $N$ satisfying $\Res^{\fg}_{\fh} N\in\cM(\fh)$, in such a way that $\Res^{\fg}_{\fh} N^{\#}= (\Res^{\fg}_{\fh} N)^{\#}$.
The map $N\iso (N^{\#})^{\#}$ is given by $v\mapsto (-1)^{p(v)}v$.
By~(\ref{eqn-frisk}),  $L(\lambda)^{\#}\cong L(\lambda)\cong \Pi L(\lambda)$ for $\lambda\in I_1$ and
\begin{equation}\label{eqn-frisk-g-mod}
L(\lambda)^{\#}\cong \Pi^{\frac{\rk F_{\lambda}}{2}} L(\lambda)\ \text{ for }
\lambda\in I_0.\end{equation}

The antiautomorphism $(-)^\#$  exists for Kac-Moody superalgebras. For $\fgl(m|n)$
the antiautomorphism $(-)^\#$  can be given by the formula
$a^{\#}:=a^t$ for $a\in\fg_{\ol{0}}$ and $a^{\#}:=\sqrt{-1}a^t$ for $a\in\fg_{\ol{1}}$
(where $a^t$ stands for the transposed matrix); this antiautomorphism on $\fg\fl(n|n)$
induces $(-)^\#$ for the
queer superalgebras $\fq_n,\fsq_n,\mathfrak{pq}_n,\mathfrak{psq}_n$.  

\subsubsection{Remark}
The duality $(-)^\#$ can be defined using a ``naive antiautomorphism'', i.e.
an invertible map $\sigma':\fg\to\fg$ satisfying $\sigma'([a,b])=[\sigma'(b),\sigma'(a)]$ via the formula $g.f(v):=f(\sigma(g)v)$ for $g\in \fg$, $f\in N^*$ and
$v\in N$ (this was done in~\cite{Gq} and~\cite{Frisk}). 

Consider the map $\theta':\fg\to\fg$ given by
$\theta'(g)=(\sqrt{-1})^mg$, where $m=0$ for $g\in \fg_{\ol{0}}$ 
and $m=1$ for $g\in \fg_{\ol{1}}$. 
If $\sigma'$ is a ``naive antiautomorphism'', then 
$\sigma'\theta'$ is an antiautomorphism.

\section{Grothendieck rings and (super)character morphisms}

Let $\CC$ be a full subcategory of the category of $\fg$-modules.
In this paper we always assume that all modules in $\CC$ are of finite length and
that $\CC$ is a dense subcategory, i.e. for every short exact sequence $0\to M'\to M\to M''\to 0$ the module $M$ lies in $\CC$ if and only if $M',M''$ lie in $\CC$.
In addition, we usually assume that $\Pi\CC=\CC$. 

We also use the following notation: for a super ring $B$ we set $B_{\mathbb{Q}}:=B\otimes_{\mathbb{Z}}\mathbb{Q}$.
\subsection{ Grothendieck groups}\label{section-groth-groups}
We denote by ${\Gr}(\CC)$ the Grothendieck group of $\CC$ which is the free abelian group generated by $N_{gr}$ for each module $N$ in $\CC$, modulo the relation that $[N]=[N']+[N'']$ whenever $0\to N'\to N\to N''\to0$ is a short exact sequence in $\mathcal{C}$.  When $\Pi\CC=\CC$, we define the structure of a $\mathbb{Z}[\xi]/(\xi^2-1)$-module on
$\Gr(\CC)$ by setting $\xi N_{gr}:=\Pi(N)_{gr}$ for $N\in\CC$.   Set $\Gr_{\pm}(\CC):=\Gr(\CC)/(\xi\mp1)$.

We call the group $\Gr_{-}(\CC):=\Gr(\CC)/(\xi+1)$ the {\em reduced Grothendieck group}.
We denote by $N_{gr,\pm}$ the image of $N_{gr}$
in $\Gr_{\pm}(\CC)$; later we will use $[N]$ for $N_{gr,-}$.

Note that $\Gr(\CC)$ is an abelian group, so $\Gr_{\pm}(\CC)$ are also abelian groups.

\subsubsection{}
We denote by $\Irr(\CC)$ the
set of isomorphism classes of irreducible modules in $\CC$ modulo $\Pi$ and
write 
$$\Irr(\CC)=\Irr(\CC)_{\ol{0}}\coprod \Irr(\CC)_{\ol{1}},$$ where $L\in\Irr(\CC)_{\ol{0}}$
if $\Pi(L)\not\cong L$ and $L\in\Irr(\CC)_{\ol{1}}$
if $\Pi(L)\cong L$.   

By our assumptions on $\CC$, $\Gr(\CC)$  a free $\mathbb{Z}$-module
with a basis 
$$\{L_{gr},\xi L_{gr}|\ L\in \Irr(\CC)_{\ol{0}}\}
\coprod \{L_{gr}|\  L\in \Irr(\CC)_{\ol{1}}\}.$$

\subsubsection{}
The  group $\Gr_+(\CC)$ is a free $\mathbb{Z}$-module
with a basis 
$\{L_{gr, +}| L\in \Irr(\CC)\}$.

\subsubsection{Reduced Grothendieck group}  One has
\begin{equation}
\Gr_{-}(\CC)=\Gr_{-}(\CC)_{free}\oplus \Gr_{-}(\CC)_{2-tor},
\end{equation}
where $\Gr_{-}(\CC)_{free}$ is a free $\mathbb{Z}$-module
with a basis 
$\{L_{gr,-} |\ L\in \Irr(\CC)_{\ol{0}}\}$ and $\Gr_{-}(\CC)_{2-tor}$ is a free $\mathbb{Z}/2\mathbb{Z}$-module
with a basis 
$\{L_{gr,-}|\ L\in \Irr(\CC)_{\ol{1}}\}$.

\begin{proposition}\label{proposition-Gr-decomp}
Consider the natural map $\psi: \Gr(\CC)\to \Gr_-(\CC)\times\Gr_+(\CC)$.
\begin{enumerate}
\item
The map $\psi$ is an embedding.
\item
The image of $\psi$ is the subgroup consisting of the pairs
$$\left(\sum_{L\in \Irr(C)} m_L L_{gr,+}, \sum_{L\in \Irr(C)} n_L L_{gr,-}\right)$$ where
$m_L,n_L\in\mathbb{Z}$ with 
$m_L\equiv n_L \mod 2$ for all $L$ and
$n_L\in \{0,1\}$ for $L\in\Irr(C)_{\ol{1}}$.
\item The  map $\psi$ induces
an isomorphism
$$\Gr(\CC)_{\mathbb{Q}}\iso \Gr_-(\CC)_{\mathbb{Q}}
\times \Gr_+(\CC)_{\mathbb{Q}}.$$
\end{enumerate}
\end{proposition}
\begin{proof}
For (i) take  $a\in\Gr(\CC)(\xi-1)\cap  \Gr(\CC)(\xi+1)$. Then
$a=a_+(\xi-1)=a_-(\xi+1)$ for $a_{\pm}\in \Gr(\CC)$. Using $\xi^2=1$ we obtain
\[
2a_+(1-\xi)=a_+(\xi-1)^2=a_-(\xi^2-1)=0.
\]
Since $\Gr_+(\CC)$ is a free
 $\mathbb{Z}$-module this gives $a_+(\xi-1)=0$, so $a=0$. This establishes (i).

The assertion (ii) follows from the observation  that  the subgroup
generated by the pairs
 $\psi(L_{gr})=(L_{gr,+}, L_{gr,-})$, $\psi(\Pi L_{gr})=(L_{gr,+}, -L_{gr,-})$
for all $L\in\Irr(\CC)$ coincides with the subgroup described in (ii).

Finally, (iii) follows from (ii).
\end{proof}

\subsubsection{Grothedieck rings}
If $\CC$ is closed under $\otimes$, then $\Gr(\CC), \Gr_{\pm}(\CC)$ 
are commutative  rings with unity
and $\psi$ in Proposition \ref{proposition-Gr-decomp} is a ring homomorphism. In this case
 $\Gr_{-}(\CC)_{2-tor}$  is an ideal in  $\Gr_-(\CC)$.

If $\CC$ is rigid, $\Gr(\CC)$ is equipped by an involution $*$ and
$\Gr_{\pm}(\CC)$ and $\Gr_-(\CC)_{free}, \Gr_-(\CC)_{2-tor}$ are $*$-invariant.

\subsection{The map $\ch_{\fh,\xi}$}\label{section-superchar-map}  Let $\fg'$ be a subalgebra of $\fg$ and let $\CC'$ be a category of $\fg'$ modules such that restriction induces a functor $\Res^{\fg}_{\fg'}:\CC\to \CC'$.   For a suitable category $\CC'$ for $\fg'$-modules, this functor induces a map $\res_{\fg'}:\Gr(\CC)\to\Gr(\CC')$ which is very useful if 
$\Gr(\CC')$ is simple enough. Below we consider this map for
the cases when $\fg'=\ft$ is a maximal toral subalgebra and for $\fg'=\fg^{\ft}=\fh$, a maximal quasitoral subalgebra.

As we will see below, $\res_{\fh}$ is an embedding if $\CC$ is ``nice enough'';
in this case  $\res_{\ft}$ induces an embedding  $\Gr_+(\CC)\to \Gr_+(\cM(\ft))$ and this map
is given by the usual (non-graded) characters.

\subsubsection{}
Let $\fh$ be quasitoral, and let $\tilde{R}(\fh)$ be the $\mathbb{Z}[\xi]$-module consisting of the sums
\[
\sum_{\nu\in I_0} (m_{\nu}+k_{\nu}\xi) [C(\nu)]+\sum_{\nu\in I_1} m_{\nu} [C(\nu)],
\ \ \ m_{\nu},k_{\nu}\in\mathbb{Z},
\]
with the $\xi$-action given by
\begin{equation}
\xi\bigl(\sum_{\nu\in I_0} (m_{\nu}+k_{\nu}\xi) [C(\nu)]+\sum_{\nu\in I_1} m_{\nu} [C(\nu)]\bigr)=\sum_{\nu\in I_0} (m_{\nu}\xi+k_{\nu}) [C(\nu)]+\sum_{\nu\in I_1} m_{\nu} [C(\nu)].\end{equation} 
For $N\in\cM(\fh)$ we introduce
\[
\ch_{\fh,\xi}(N):=\sum_{\nu\in I_0} (m_{\nu}+k_{\nu}\xi) [C(\nu)]+\sum_{\nu\in I_1} m_{\nu} [C(\nu)]\in \tilde{R}(\fh)
\]
where $m_{\nu}:=[N_{\nu}:C(\nu)]$ and $k_{\nu}:=[N_{\nu}:\Pi(C(\nu))]$.  

This defines a linear map $\ch_{\fh,\xi}:\Gr(\cM(\fh))\to \tilde{R}(\fh)$, which we refer to as the graded $\fh$-character of $N$. We denote $\sch N$ the image of $\ch_{\fh,\xi}$ in
$\tilde{R}(\fh)/\tilde{R}(\fh)(\xi+1)$. Then
$$ \sch_{\fh}(N):=\sum_{\nu\in I_0} (m_{\nu}-k_{\nu}) [C(\nu)]+\sum_{\nu\in I_1} m_{\nu}' [C(\nu)]$$
where $m'_{\nu}=0$ if $m_{\nu}$ is even and $m'_{\nu}=1$ if $m_{\nu}$ is odd.

\begin{lemma}\label{lemma-realization-Gr-h} The maps $[N]\to \ch_{\fh,\xi} N$, respectively $[N]\to\sch_{\fh} N$ define isomorphisms
$\Gr(\cM(\fh))\to \tilde{R}(\fh)$, respectively  $\Gr_-(\cM(\fh))\to \tilde{R}(\fh)/\tilde{R}(\fh)(\xi+1)$. Further, these maps are compatible
with $(-)^{\#}$. 
\end{lemma}

\subsubsection{Note:} Because of this lemma, we will subsequently do away with the notation $\tilde{R}(\fh)$ and instead directly identify $\Gr(\cM(\fh))$ (and $\Gr_-(\cM(\fh))$) with the corresponding spaces $\tilde{R}(\fh)$ (and $\tilde{R}(\fh)/\tilde{R}(\fh)(1+\xi)$) as presented above.

\begin{proof}[Proof of Lemma \ref{lemma-realization-Gr-h}]
	We show that $\ch_{\fh,\xi}$ is an isomorphism, with the result for $\sch_{\fh}$ following easily.  Clearly $\ch_{\fh,\xi}$ is surjective, so it suffices to prove injectivity.  First we observe that for $N$ in $\cM(\fh)$, we have the following equality in the Grothendieck ring $\Gr(\cM(\fh))$:
	\[
	[N]=\left[\bigoplus\limits_{\nu\in I_0}\left(C(\nu)^{\oplus m_{\nu}}\oplus \Pi C(\nu)^{\oplus k_{\nu}}\right)\oplus\bigoplus\limits_{\nu\in I_1} C(\nu)^{m_{\nu}}\right]
	\]
	This simply follows from the fact that $N$ has finite Loewy length, since this is true for the algebras $\Cl(\lambda)$.  Further, it is not difficult to see that a basis of $\Gr(\cM(\fh))$ is given by elements of the form:
\[
	\left[\bigoplus\limits_{\nu\in I_0}\left(C(\nu)^{\oplus m_{\nu}}\oplus \Pi C(\nu)^{\oplus k_{\nu}}\right)\oplus\bigoplus\limits_{\nu\in I_1} C(\nu)^{m_{\nu}}\right]
\]
where $m_{\nu},k_{\nu}\in\mathbb{N}$.  From this the isomorphism easily follows.  The compatibility with $(-)^{\#}$ is obvious.  
\end{proof}

\begin{corollary}  
\begin{enumerate}
\item
One has $\ch_{\fh,\xi} \Pi N=\xi\ch_{\fh,\xi} N$, $\ \ \sch_{\fh} \Pi N=-\sch_{\fh} N$.

\item For $\lambda\in I_1$ one has
$\sch_{\fh} L(\lambda)=\sum_{\mu\in I_1} m_{\mu} [C(\mu)]$.
\item
If $\fg$ admits $(-)^{\#}$ as in Section \ref{section-duality}, then
the coefficient of $[C(\nu)]$ in
$\sch_{\fh} L(\lambda)$ is zero if $\rk F_{\lambda},\rk F_{\nu}$ are even and
$\rk F_{\nu}\not\equiv \rk F_{\lambda} \mod 4$.
\end{enumerate}
\end{corollary}

\begin{proof}
The assertions follow from~(\ref{eqn-I01})
and~(\ref{eqn-frisk-g-mod}). 
\end{proof}

The next corollary is a direct generalization of~\cite{SV}, Proposition 4.2.

\begin{corollary}\label{corollary-res-Gr-embedding}
Let $\CC$ be a full subcategory of $\fg$-modules
with the following properties: each module in $\CC$ is of 
 finite length and is locally finite over $\fb$, and the restriction to $\fh$ lies in $\cM(\fh)$.  Then the map $\Res^{\fg}_{\fh}$ induces injective maps
$\ch_{\fh,\xi}:\Gr(\CC)\hookrightarrow \Gr(\cM(\fh))$ and $\sch_{\fh}:\Gr_-(\CC)\hookrightarrow \Gr_-(\cM(\fh))$.

\end{corollary}

\begin{proof} 
Let us check the injectivity of the first map
$\ch_{\fh,\xi}: \Gr(\CC)\to \Gr(\cM(\fh))$.
Any simple module in $\CC$ is  $L(\lambda)$ or  $\Pi(L(\lambda))$ for some $\lambda\in\ft^*$. Since every module in $\CC$ has finite length,
$\Gr(\CC)$ is a free $\mathbb{Z}$-module spanned by $[L(\lambda)],\xi [L(\lambda)]$ for $\lambda\in I_0$ and $[L(\lambda)]$ for $\lambda\in I_1$.  Assume that
\[
\ch_{\fh,\xi}\bigl(\sum_{i=1}^s (m_i+k_i\xi) [L(\lambda_i)]\bigr)=0
\]
where $k_i=0$ for $\lambda_i\in I_1$ and
$\gamma(\lambda_1)$ is maximal among $\gamma(\lambda_i)$ for $i=1,\ldots,s$.
Then for $i=2,\ldots, s$ one has
 $L(\lambda_i)_{\lambda_1}=0$, so  the coefficient of $[C(\lambda_1)]$
in $\ch_{\ft}([L(\lambda_i)]$ is zero. Hence
$(m_1+k_1\xi)[L(\lambda_1)]=0$. This gives $m_1=k_1=0$
and implies the injectivity of $\ch_{\fh,\xi}$. The injectivity
of $\sch_{\fh}$ easily follows. \end{proof}

\subsection{Example}
The category of  finite-dimensional $\fg$-modules $\cF in(\fg)$ is a  rigid tensor category with the duality $N\mapsto N^*$
given by the anti-automorphism $-\Id_{\fg}$. 
Note that the Grothendieck ring $\Gr({\cF}in(\fg))$
 is a commutative ring with a basis $\{[L]\}$, where
$L$ runs through isomorphism classes of finite-dimensional simple modules. Corollary \ref{corollary-res-Gr-embedding} gives an embedding $\ch_{\fh,\xi}: \Gr(\cF in(\fg))\ \xhookrightarrow\ \Gr(\cF in(\fh)))$ ($[M]\mapsto \ch_{\fh,\xi} M$) which is a ring homomorphism.  By abuse of notation we will also denote the image of this homomorphism by $\Gr(\cF in(\fg))$. 

The duality induces an involution on $\Gr({\cF}in(\fg))$, which we also denote by $*$.
One has $\xi^*=\xi$. The homomorphism  $\ch_{\fh,\xi}: \Gr({\cF}in(\fg))\hookrightarrow\ \Gr(\cF in(\fh))$ is compatible with $*$, so $\Gr(\cF in(\fg))$ is a $*$-stable subring of $\Gr(\cF in(\fh))$.

\subsubsection{Remark}
Let $\fg$ be a Kac-Moody superalgebra (so $\ft=\fh$) and
let $\Lambda_{int}\subset\fh^*$ be a lattice contatining $\Delta(\fg)$ such that
the parity $p: \Delta\to\mathbb{Z}_2=\{\ol{0},\ol{1}\}$ can be extended to
$p:\Lambda_{int}\to\mathbb{Z}_2$.	Assume, in addition, that for the category $\CC$, each $N\in \CC$ has $N_{\nu}=0$ if $\nu\not\in\Lambda_{int}$.

Then $\CC=\CC_{+}\oplus \Pi(\CC_{+})$, where $N\in\CC$ lies in
$\CC_{+}$ if and only if $N_{\nu}\subset N_{p(\nu)}$.  We have that
\[
\Gr(\CC)=\Gr(\CC_{+})\times \mathbb{Z}[\xi]/(\xi^2-1).
\]
and thus one can recover $\Gr(\CC)$ from $\Gr(\CC_+)$; however we note that $\CC_+$ is not $\Pi$-invariant.

Further, we have in this case that $\Gr_-(\CC)\cong \Gr_+(\CC)\cong \Gr(\CC_+)$.  If $\CC$ is a tensor category, then $\CC_+$ is also a tensor category
(but $\Pi(\CC_{+})$ is not).  In~\cite{SV}, Sergeev and Veselov described the ring $\Gr(\CC_+)$
for the finite-dimensional Kac-Moody superalgebras.

\subsection{The ring $\Gr(\fh)$}  Let $\fh$ be quasitoral.  We write $\Gr(\fh):=\Gr(\cF in(\fh))$.  
The map $[N]\mapsto \ch_{\fh,\xi} N$ introduced in Section \ref{section-superchar-map} gives an isomorphism
of  $\Gr(\fh)$ and the free $\mathbb{Z}$-module
spanned by $[C(\nu)]$, $\xi[C(\nu)]$ for $i\in I_0$ and
$[C(\nu)]$ for $i\in I_1$. 
We view $\Gr(\fh)$ as a commutative algebra endowed with the involutions
$(-)*$ and $(-)^{\#}$. One has $C(\lambda)\in \Irr(\CC)_{\ol{i}}$ if and only if
$\lambda\in I_i$.

One has $[E_0]=1$, $[\Pi(E_0)]=\xi$ and 
\renewcommand{\arraystretch}{1.6}
\begin{equation}
\begin{array}{l}
\ [C(\lambda)]^* \in \{ [C(-\lambda)], \xi [C(-\lambda)]\},\ \ \
 [C(\lambda)]^{\#}=\xi^{\frac{\rank{F_{\lambda}}}{2} }[C(\lambda)], 
\text{ for }\lambda\in I_0\\
\ \xi [C(\lambda)]=[C(\lambda)], \ \ \ [C(\lambda)]^*=[C(-\lambda)], \ \ \ 
 [C(\lambda)]^{\#}=[C(\lambda)], 
\text{ for }\lambda\in I_1.\end{array}
\end{equation}

\subsubsection{}
The multiplication in $\Gr_+(\fh)$ is given by
\[
C(\lambda)_{gr,+}C(\nu)_{gr,+}=\frac{\dim C(\lambda)\dim C(\nu)}{\dim C(\lambda+\nu)}C(\lambda+\nu)_{gr,+}.
\]

Let $\mathbb{Z}[e^{\nu},\nu\in\ft^*]$ be the group ring of $\ft^*$.
For $N\in\cM(\fh)$ we set
$$\ch_{\ft} N:=\sum_{\nu\in \ft^*}\dim N_{\nu} e^{\nu}.$$ 
The map $N\mapsto \ch_{\ft} N$ induces an embedding  
$\Gr_+(\fh)\hookrightarrow\mathbb{Z}[e^{\nu},\nu\in\ft^*]$.
The image is the subring of elements $\sum_{\nu} m_{\nu}e^{\nu}$ where 
$m_{\nu}$ is divisible by $\dim C(\nu)$, and we have an isomorphism 
$\Gr_+(\fh)_{\mathbb{Q}}\iso \mathbb{Q}[e^{\nu},\nu\in\ft^*]$.

\subsection{Spoiled superalgebras}\label{section-spoiled-algs}  We call a $\mathbb{Z}$-superalgebra $A=A_{\ol{0}}\oplus A_{\ol{1}}$ {\em spoiled} if $2A_{\ol{1}}=A_{\ol{1}}^2=0$.  Given a superalgebra $B$, we make it into a spoiled superalgebra $B^{spoil}$ by setting, as a ring:
\[
B^{spoil}=\left(B\otimes_{\mathbb{Z}}\mathbb{Z}[\vareps]/(\vareps^2,2\vareps)\right)_{\ol{0}}.
\]
Here $\vareps$ is an element of degree 1.  The $\mathbb{Z}_2$-grading on $B^{spoil}$ is declared to be $B^{spoil}_{\ol{0}}=B_{\ol{0}}$ and $B^{spoil}_{\ol{1}}=B_{\ol{1}}\vareps$.  We in particular observe that 
\[
B^{spoil}_{\mathbb{Q}}:=(B_{\ol{0}})_{\mathbb{Q}}.
\]

\subsubsection{The algebra $\Gr_-(\fh)$}\label{section-reduced-Gr-h}  
From now on we will write $[M]$ in place of $M_{gr,-}$.  Using Corollary \ref{corollary-tensor-prod-parity-invce} we obtain that in $\Gr_-(\fh)$ we have the following multiplication law:
$$[C(\nu)][C(\lambda)]=\left\{\begin{array}{ll}
\pm [C(\lambda+\nu)] &  {\text{ if }
{\rank F_{\lambda}+\rank F_{\nu}=\rank F_{\lambda+\nu}, \
 \rank F_{\lambda}\cdot\rank F_{\nu}\equiv 0\mod 2}} \\
0 &  \text{ otherwise}.\end{array}\right.$$

As a result,  $\Gr_-(\fh)$ is $\mathbb{Z}$-graded algebra
\begin{equation}\label{eqn-grading-red-Gr-h}
\Gr_-(\fh)=\bigoplus_{i=0}^{\infty} \Gr_-(\fh)_i
\end{equation}
 where $\Gr_-(\fh)_i$ is spanned by $[C(\nu)]$ with 
$\rank F_{\nu}=i$.   We consider the corresponding
$\mathbb{Z}_2$-grading
$$\Gr_-(\fh)_{\ol{0}}:=\bigoplus_{i=0}^{\infty} \Gr_-(\fh)_{2i},\ \ \ 
\ \Gr_-(\fh)_{\ol{1}}:=\bigoplus_{i=0}^{\infty} \Gr_-(\fh)_{2i+1}.$$
One has $\Gr_-(\fh)_{\ol{1}}\cdot \Gr_-(\fh)_{\ol{1}}=0$ and $2\Gr_-(\fh)_{\ol{1}}=0$, so that $\Gr_-(\fh)$ is a spoiled algebra.  The following corollary is clear.

\begin{corollary}
	\begin{enumerate}
		\item The algebra $\Gr_-(\fh)$ is a spoiled superalgebra with $$\begin{array}{l}\Gr_-(\fh)_{\ol{0}}=\bigoplus_{i=0}^{\infty} \Gr_-(\fh)_{2i}=\Gr_-(\fh)_{free}\\ \Gr_-(\fh)_{\ol{1}}=\bigoplus_{i=0}^{\infty} \Gr_-(\fh)_{2i+1}=\Gr_-(\fh)_{2-tor}.\end{array}$$ 		
		\item
		The algebra $\Gr(\fh)$ is isomorphic to a subalgebra of
		$\mathbb{Z}[e^{\nu},\nu\in\ft^*]\times \Gr_-(\fh)$
		consisting of
		$$(\sum_{\nu\in\ft^*} m^+_{\nu} \dim C(\nu) e^{\nu}, \sum_{\nu\in I_0} m^-_{\nu} [C(\nu)]+ \vareps\sum_{\nu\in I_1} m^-_{\nu}[C(\nu)])$$
		where $m^{\pm}_{\nu}\in\mathbb{Z}$ are equal to zero except for finitely many
		values of $\nu$, $m^+_{\nu}\equiv m^-_{\nu}$ modulo $2$
		and $m^-_{\nu}\in\{0,1\}$ for $\nu\in I_1$.
		
		\item The algebra $\Gr(\fh)_{\mathbb{Q}}$ is isomorphic to $\mathbb{Q}[e^{\nu},\nu\in\ft^*]\times \Gr_-(\fh)_{\mathbb{Q}}$.
	\end{enumerate}
\end{corollary}

\subsubsection{}
 For the rest of this section we set $\Gr(\fg):=\Gr(\cF in(\fg))$, where $\fg$ is as in Section 4.  Let
 $\psi_{\pm}(\fg): \Gr(\fg)\to\Gr_{\pm}(\fg)$ be the canonical 
epimorphisms.
By Proposition \ref{proposition-Gr-decomp}, $\psi_+\times\psi_-$ gives an embedding
$\Gr(\fg)\hookrightarrow \Gr_+(\fg)\times \Gr_-(\fg)$.

We will use the following construction:
for any subsets  $A_{\pm}\subset \Gr_{\pm}(\fh)$ we introduce
$$A_+\underset{\Gr(\fh)}{\times} A_-:=\{a\in \Gr(\fh)|\ \psi_{\pm}(a)\in A_{\pm}\}.$$
Note that  $A_+\underset{\Gr(\fh)}{\times} A_-$ 
 is a subring of $\Gr(\fh)$ if $A_{\pm}$ are rings.
 
\begin{lemma}
Let $A$ be a $\mathbb{Z}[\xi]$-submodule of $\Gr(\fh)$
with the following property: if $a\in \Gr(\fh)$ and $2a\in A$, then $a\in A$.
Then
$$A=\psi_+(A)\underset{\Gr(\fh)}{\times} \psi_-(A).$$
\end{lemma}

\begin{proof}
Take 
$a\in \Gr(\fh)$ such that $\psi_{\pm}(a)\in \psi_{\pm}(A)$. Then $A$ contains
$a-c(1-\xi)$  for some $c\in \Gr(\fh)$.
Since $A$ is a  $\mathbb{Z}[\xi]$-submodule, $A$ contains 
$(1+\xi)(a-c(1-\xi))=(1+\xi)a$. Similarly, $A$ contains $(1-\xi)a$, so
 $2a\in A$. Then the assumption gives $a\in A$ as required.
\end{proof}

\begin{corollary}
$\ \Gr(\fg)=\Gr_+(\fg)\underset{\Gr(\fh)}{\times}\Gr_-(\fg)$.
\end{corollary}

\begin{proof}
Recall by Corollary \ref{corollary-res-Gr-embedding} that $\Gr(\fg)$ is a $\mathbb{Z}[\xi]$-subring of $\Gr(\fh)$.
Take
 $a\in \Gr(\fh)$ with $2a\in \Gr(\fg)$.
Write 
	$$2a=\sum_{\lambda\in P^+(\fg)} m_{\lambda} \ch_{\fh,\xi} L(\lambda)
=\sum_{\nu } 2m'_{\nu}[C(\nu)]$$
where $m_{\lambda}, m'_{\nu}\in \mathbb{Z}[\xi]$, and $P^+(\fg)$ is the set of dominant weights of $\fg$.  Let $\lambda$
be maximal such that $m_{\lambda}\not=0$. Then $m_{\lambda}=2m'_{\lambda}$, so we may subtract $2L(\lambda)_{gr}$, and conclude by induction.
\end{proof}

\renewcommand{\arraystretch}{1}

\subsection{Equivariant setting}\label{section-equivariant-setting}  Suppose that $G$ is a finite group which acts on a quasitoral superalgebra $\fh$ by automorphisms.  Our main example of this setup is when we consider quasireductive Lie superalgebras in Section \ref{section-quasired}, and $G$ is the Weyl group.

The group $G$ then acts naturally on $\Gr(\fh)$ by twisting, i.e. $g\cdot V_{gr}=V_{gr}^g$ for an $\fh$-module $V$.  This descends to a natural action of $G$ on $\Gr_-(\fh)$.  Suppose that $\nu\in\ft^*$ and $g\in G$ such that $g\nu=\nu$.  Then we have either $g\cdot C(\nu)_{gr}=C(\nu)_{gr}$ or $g\cdot C(\nu)_{gr}=\xi C(\nu)_{gr}$.  It follows that on the reduced Grothendieck ring we have $g\cdot [C(\nu)]=\pm[C(\nu)]$. Thus $\Stab_G[C(\nu)]\sub\Stab_G\nu$.  Define
\[
\ft_{G}^*=\{\nu\in\ft^*:\Stab_G[C(\nu)]=\Stab_G\nu\}.
\]
We observe that $\ft_{G}^*$ is a non-empty, $G$-stable cone in $\ft^*$.  It need not be open or closed, and it may consist only of $0$. It is clear that $I_1\sub\ft_G^*$. 

Let $\nu\in \ft^*$.  We introduce the grading on $C(\nu)$ in the  way 
described in Remark \ref{remark-rigidify-parity}.  Since $\ol{T}_{\nu}$ is proportional to a product of a basis elements in $\fh_{\ol{1}}/K_{\lambda}$, one has 
\begin{equation}
g(\ol{T}_{\nu})=\det (g|_{\fh_{\ol{1}}/K_{\nu}}) \ol{T}_{\nu}\ \ \text{ 
	for each }g\in \Stab_W\nu.
\end{equation} 
Since $g$ is acting by an orthogonal transformation on $\fh_{\ol{1}}/K_{\nu}$ we have $\det(g|_{\fh_{\ol{1}}/K_{\nu}})=\pm1$.  Therefore for $\nu\in I_0$ and $g\in \Stab_G\nu$ we have
\[
g[C(\nu)]=\det (g|_{\fh_{\ol{1}}/K_{\nu}}) [C(\nu)].
\]

\begin{corollary}\label{corollary-equivar-red-Gr-struc} 
		\begin{enumerate}
		\item We have:
		\[
		\ft_{G}^*\cap I_0=\{\nu\in I_0:\det(g|_{\fh_{\ol{1}}/K_{\nu}})=1\text{ for all }g\in\Stab_{G}\nu\}.
		\]		
		\item For $\nu\in\ft_G^*$, the element  
		\[
		a_{\nu}:=\sum_{g\in G/\Stab_G \nu} g[C(\nu)]
		\]
		is non-zero and well-defined.  
		\item 	The algebra $\Gr_-(\fh)^G$ is naturally a spoiled superalgebra; the even part has a $\mathbb{Z}$-basis given by $a_{\nu}$ for a choice of coset representatives $\nu\in \ft_{G}^*\cap I_0/G$.  The odd part has a $\mathbb{Z}_2$-basis given by $a_{\nu}$ for a choice of coset representatives $\nu\in I_1/G$.
	\end{enumerate}
\end{corollary}

\section{$\fh$-supercharacters of some highest weight $\fg$-modules}
We continue to write $[N]$ for the image a module $N$ in the corresponding
reduced Grothendieck group.

Fix a triangular decomposition $\fg=\fn^-\oplus\fh\oplus\fn^+$ coming from $\gamma:\mathfrak{t}^*\to\mathbb{R}$ as in Section 4, and consider the corresponding category $\CO$ with respect to $\fb=\fh\oplus\fn^+$.  To be precise, $\CO$ consists of all finitely generated $\fg$-modules which are weight modules and $\fn$-locally finite.  We take $M(\lambda)$ as in Section \ref{section-highest-wt-mod}.  For weights $\lambda,\mu\in\ft^*$ we write $\lambda\sim_0\mu$ if $L(\lambda)$ is a composition factor of $M(\mu)$; then we let $\sim$ be the equivalence relation on $\mathfrak{t}^*$ generated by $\sim_0$.

The goal of this section is to prove Theorem \ref{theorem-parity_invce}, which in some sense is the best version of Theorem \ref{theorem-intro_1} that holds in great generality.  The idea is to enforce assumptions that guarantee that all Verma modules with highest weights lying in a fixed equivalence class of interest have an especially nice $\fh$-supercharacter.  

\subsection{Notation} Let $\fa'\sub\fa$ be Lie superalgebras, $L'$ a simple $\fa'$-module, and $N$ an $\fa$-module such that $\Res^{\fa}_{\fa'}N$ has a finite length.  Set
\begin{equation}
\smult(N:L'):=\left\{\begin{array}{ll}
\ [\Res^{\fa}_{\fa'} N:L']-[\Res^{\fa}_{\fa'} N:\Pi L'] & \text{ if }L'\not\cong\Pi L' \\
\ [\Res^{\fa}_{\fa'} N:L'] \text{ mod } 2 & \text{ if }L'\cong\Pi L'.\end{array}\right.
\end{equation}

If $\Res^{\fg}_{\fh} N\in\cM(\fh)$ we set
\begin{equation}
\Omega(N):=\{\mu\in\ft^*| N_{\mu}\not=0\},\ \ \ \ 
\sOmega(N):=\{\mu\in\ft^*|\smult(N_{\mu}:C(\mu))\neq0 \}.
\end{equation}
Then
$$\sch_{\fh} N=\sum_{\nu\in \sOmega(N)} \smult(N_{\nu}:C(\nu))[C(\nu)].$$

\subsection{On $\sch_{\fh} M(\lambda)$}
For $\nu\in Q^-$ we have an isomorphism of $\fh$-modules:
	\[
	M(\lambda)_{\lambda+\nu}\cong\mathcal{U}(\fn^-)_{\nu}\otimes C(\lambda).
	\]

By Theorem \ref{theorem-tensor-prod}	, $\sch_{\fh} (C(\nu)\otimes C(\lambda))\not=0$ implies
$\dim C(\nu)\cdot \dim C(\lambda)=\dim C(\lambda+\nu)$. Therefore
\begin{equation}\label{eqn-supersupport_verma}
\sOmega (M(\lambda))\subset \{\nu\in\ft^*|\ 
\dim C(\nu)\cdot \dim C(\lambda)=\dim C(\lambda+\nu)\}.\end{equation}

\begin{corollary}
Assume that $F_{\nu}\not=0$ for $\nu\in Q^-\setminus\{0\}$.
If  $\corank F_{\lambda}\leq 1 $ and $\dim\fh_{\ol{1}}$ is even, or $\corank F_{\lambda}=0$, then
$$\sch_{\fh} M(\lambda)=\sch_{\fh} C(\lambda).$$
\end{corollary}

\begin{proof}
If $\corank F_{\lambda}=0$, or $\corank F_{\lambda}=1$ and $\dim\fh_{\ol{1}}$ is even, then 
$\dim C(\lambda)\geq \dim C(\lambda+\nu)$ for any $\nu$, so the
 formula follows from~(\ref{eqn-supersupport_verma}).
\end{proof}

\subsection{$\sch_{\fh}L(\lambda)$ when $\corank F_{\lambda}\leq 1$}  

We make the following assumptions on our triangular decomposition:
\begin{enumerate}
	\item[(A1)] Setting $Q^-=\mathbb{N}\Delta^-$, we have that $\gamma(Q^-)\sub\mathbb{R}^-$ is discrete.
	
	\item[(A2)]  One has $F_{\nu}\neq 0$ for each
$\nu\in Q^{-}\setminus\{0\}$.
		
	\item[(A3)] For $\lambda\in\ft^*$ such that $\operatorname{corank} F_{\lambda}\leq1$, we have that $M(\lambda)$ has finite length as a $\fg$-module, (one can weaken this to assume that it admits a ``local composition series'' in the sense of~\cite{DGK}).
	\end{enumerate}

\begin{theorem}\label{theorem-parity_invce}
Let $\Lambda\subset \ft^*$ be an equivalence class of $\sim$ such that
\begin{enumerate}
	\item[(A4)]
$\corank F_{\lambda}\leq 1$ for each $\lambda\in\Lambda$ if $\dim\fh_{\ol{1}}$ is even, and otherwise $\corank F_{\lambda}=0$ for each $\lambda\in\Lambda$. 
	\end{enumerate}
Then
for each $\lambda\in\Lambda$ one has
\begin{equation}
\sch_{\fh} L(\lambda)=\sum_{\mu\in\Lambda} m_{\mu} \sch_{\fh} C(\mu)
\end{equation}
for some $m_{\mu}\in\mathbb{Z}$.
\end{theorem}

\begin{proof}
Assume that the assertion does not hold. Take $(\mu,\mu')$ such that 
$\gamma(\mu'-\mu)$ is maximal  with the properties
$$\mu\in \Lambda,\ \  \mu'\not\in \Lambda,\ \ \sch_{\fh} L(\mu)_{\mu'}\not=0$$
(the existence of such pair follows from (A1)).

By (A3), in the Grothendieck ring of $\fg$ we may write:
$$
[L(\mu)]=[M(\mu)]-\sum_{\nu\in Q^-\setminus\{0\}: \mu+\nu\in \Lambda}
a_{\mu\nu}[L(\mu+\nu)]$$
where  $a_{\mu\nu}\in\mathbb{N}$.  We may map the above formula to the Grothendieck ring of $\fh$-modules via restriction, and thus learn that 
$$\sch_{\fh} L(\mu)_{\mu'}=\sch_{\fh} M(\mu)_{\mu'}-
\sum_{\nu\in Q^-\setminus\{0\}: \mu+\nu\in \Lambda}\sch_{\fh} L(\mu+\nu)_{\mu'}.$$
By~(\ref{eqn-supersupport_verma}), $\sch_{\fh} M(\mu)_{\mu'}=0$ and the maximality of
$\gamma(\mu'-\mu)$ implies $\sch_{\fh} L(\mu+\nu)_{\mu'}=0$.
Hence  $\sch_{\fh} L(\mu)_{\mu'}=0$, a contradiction.
\end{proof}

We will see in Section \ref{section-red-Gr-q} that Theorem \ref{theorem-parity_invce} holds for any irreducible highest weight representation of $\fq_n$, without restrictions on $\operatorname{corank}F_{\lambda}$.

\section{On $\Gr_-(\fg)$ in the case when $\fg$ is quasireductive}
In this section we assume $\fg$ is quasireductive, i.e.~$\fg$ is finite-dimensional,
$\fg_{\ol{0}}$ is reductive, and $\fg_{\ol{1}}$ is a semisimple $\fg_{\ol{0}}$-module
(see~\cite{Sqred}, \cite{M}, and \cite{IO} for examples and a partial classification of such algebras). The maximal toral subalgebras $\ft$ in $\fg_{\ol{0}}$ are the Cartan
subalgebras in $\fg_{\ol{0}}$, and the maximal quasitoral subalgebras $\fh$ in $\fg$ are the Cartan
subalgebras in $\fg$; all such subalgebras are conjugate to one another by inner automorphisms because the same is true on the even part (see \cite{H}).
We fix such $\ft$ and $\fh$. Let $P(\fg)\subset \ft^*$ be the set of weights
appearing in finite-dimensional $\fg$-modules, which is the same as the set of weights appearing in finite-dimensional $\fg_{\ol{0}}$-modules.  We fix a positive system of roots on $\fg$, and thus also on $\fg_{\ol{0}}$. We denote by $P^+(\fg)$ the set of dominant weights, i.e.
\[
P^+(\fg):=\{\lambda\in\ft^*|\ \dim L(\lambda)<\infty\}.
\]	
We have
\[
P(\fg)=P^+(\fg)+\mathbb{Z}\Delta=P^+(\fg_{\ol{0}})+\mathbb{Z}\Delta_{\ol{0}}.
\]
Write $\CC$ for the subcategory of $\cF(\fh)$ consisting of modules with weights lying in $P(\fg)$. 

\subsection{On $\sch_{\fh}(\cF(\fg))$}\label{section-quasired} Define $\cF(\fg)$ to be the full subcategory of $\cF in(\fg)$ consisting of modules which are semisimple over $\fg_{\ol{0}}$.  

The Weyl group $W$ of $\fg_{\ol{0}}$ acts naturally on $\fh$, and thus we are in the setup of Section \ref{section-equivariant-setting}; we refer to that section for the definition of $\ft_{W}^*$.  We see that $W$ preserves the subcategory $\CC$.  We set
\[
P(\fg)'=\ft_{W}^*\cap P(\fg), \ \ \ \ P^+(\fg_{\ol{0}})'=\ft_{W}^*\cap P^+(\fg_{\ol{0}}).
\]
Recall that 
\[
P(\fg)'\cap I_0:=\{\lambda \in P(\fg)\cap I_0|\ \forall w\in \Stab_W \lambda, \ \ \det (w|_{\fh_{\ol{1}}/\Ker F_{\lambda}})=1\}.
\]
From the theory of reductive Lie algebras, in this case we have a natural bijection $P^+(\fg_{\ol{0}})'\to P(\fg)'/W$.  Recall that for $\nu\in\ft_W^*$ we set
\[
a_{\nu}=\sum\limits_{w\in W/\Stab_W\nu}w[C(\nu)].
\] 

\begin{theorem}\label{theorem-qred-sch}
\begin{enumerate}
\item 
	The algebra $\Gr_-(\CC)^W$ is naturally a spoiled superalgebra; its even part has $\mathbb{Z}$-basis given by $a_{\nu}$ for $\nu\in P^+(\fg_{\ol{0}})'\cap I_0$
	and its odd part has a  $\mathbb{Z}_2$-basis  $a_{\nu}$ for $\nu\in\left(P^+(\fg_{\ol{0}})\cap I_1\right)$.
\item $\sch_{\fh}$ defines an embedding $\Gr_-(\fg)\hookrightarrow \Gr_-(\CC)^W$.
\item For $\lambda\in P^+(\fg)\cap I_0$ one has 
\[
\sch_{\fh} L(\lambda)=\sum_{\nu} k_{\nu} a_{\nu}
\]
where $k_{\nu}\in\mathbb{Z}$ with  $k_{\lambda}=1$.
\item If $\lambda\in P^+(\fg)\cap I_1$, then
\[
\sch_{\fh} L(\lambda)=\sum_{\nu\in I_1} k_{\nu} a_{\nu}.
\]
with $k_{\lambda}=1$.
\item
$P^+(\fg)\subset P^+(\fg_{\ol{0}})'$.
\end{enumerate}
\end{theorem}

\begin{proof}
	Part (i) follows from Corollary \ref{corollary-equivar-red-Gr-struc}.  For part (ii) let $L:=L(\lambda)$ be a simple finite-dimensional $\fg$-module.   Clearly all its weights lie in $P(\fg)$.  Because $W$ can be realized from inner automorphisms of $\fg$, it is clear that the $\sch_{\fh}L(\lambda)$ must be $W$-invariant, proving (ii).  
	
	Part (iii) follows uses (\ref{eqn-frisk-g-mod}), and part (4) uses that $L(\lambda)\cong\Pi L(\lambda)$.  Finally part (v) then follows from (iii).
\end{proof}

\subsubsection{Remark}  The algebra $\Gr_{-}(\fg)$ is not naturally spoiled with respect to the embedding into $\Gr_-(\fh)$.  For example if $\fg=\fq_2\times\fq_1$, then consider the module $S^{2}V_{std}^{(1)}\boxtimes V_{std}^{(2)}$, where $V_{std}^{(1)},V^{(2)}_{std}$ denote the standard modules of $\fq_2,\fq_1$ respectively.  Then this module is of the form $L\oplus \Pi L$ for some simple $\fq_2\times\fq_1$-module $L$.  One may check that $\sch_{\fh}L$ is not a homogeneous element of $\Gr_-(\fh)$, with respect to its spoiled grading (see Section \ref{section-reduced-Gr-h}).

\begin{corollary}\label{corollary-longest-elt-det}
Let $\lambda\in P^+(\fg)\cap I_0$ and $w_0\in W$ be such that $w_0\lambda=-\lambda$ and $L(\lambda)^*$ has highest weight $\lambda$.  Then we have
$$L(\lambda)^*\cong \Pi^i L(\lambda)\ \  \text{ where }\ \ 
(-1)^i
=\det (w_0|_{\fh_{1}/\Ker F_{\lambda}}).$$
\end{corollary}
\begin{proof} Recall that $\phi_{\lambda}: \cU(\fh)\to \Cl(\lambda)$, $p_{\lambda}:\Cl(\lambda)\to\Cl(\lambda)/K_{\lambda}$ stand for the canonical
epimorphisms.  By Section \ref{section-dualities-fd} we have
$\ol{T}_{\lambda}=H_{2k}'\ldots H'_1$, where
$H'_1,\ldots,H_{2k}'$ is a basis
of  $p_{\lambda}\phi_{\lambda}(\fh_{1})\subset \Cl(\lambda)/K_{\lambda}$ 
satisfying $[H'_i,H'_j]=0$ for $i\not=j$.
Set
\[
\ol{T}_{-\lambda}=\ol{T}_{w_0\lambda}=w_0({H}'_{2k})\ldots w_0({H}'_1).
\]
From this definition we have $C(-\lambda)=C(\lambda)^{w_0}$.  We have
$\ol{\sigma_{-\lambda}}(H'_i)=-H'_i$, and thus
\[
\ol{\sigma_{-\lambda}}(\ol{T}_{-\lambda})=(-1)^{k(2k-1)}
w_0({H}'_1)\ldots w_0({H}'_k)=(-1)^{k(2k-1)}\det(w_0|_{\fh_{\ol{1}}/K_{\lambda}})
H_1'\ldots H'_{2k}=(-1)^i \ol{T}_{\lambda}.
\]
It follows that $C(\lambda)^*=\Pi^iC(-\lambda)$.   Therefore
\[
(L(\lambda)^*)_{\lambda}=(L(\lambda)_{-\lambda})^*=(\Pi^iC(\lambda)^*)^*=\Pi^iC(\lambda).
\]
and we obtain $L(\lambda)^*\cong \Pi^i L(\lambda)$ as required.
\end{proof}

\section{The $\DS$-functor and the reduced Grothendieck group}

In Sections \ref{section-DS-functor} and \ref{section-ds-on-red-groth}, $\fg$ is any Lie superalgebra. 
We fix $x\in\fg_{\ol{1}}$ with
$[x,x]=c\in\fg_{\ol{0}}$ such that $\ad c$ is semisimple (such elements $x$ are called \emph{homological}).

\subsection{$\DS$-functors: construction and basic properties}\label{section-DS-functor}
The $\DS$-functors were introduced in \cite{DS}; we use a slight generalization (see \cite{EAS} for a more in-depth treatment).
For a $\fg$-module $M$ and $u\in\fg$ we set $$M^u:=\Ker_M u.$$  

Let $M$ be a $\fg$-module on which $c$ acts semisimply.  Write $DS_xM=M_x:=M^x/(\im x\cap M^x)$.  Then $\fg^x$ and $\fg_x$ are Lie superalgebras, where $x$ acts via the adjoint action.

Observe that $M^x, xM^c$ are $\fg^{x}$-invariant and $[x,\fg^c] M^x\subset xM^c$,
so $\DS_x(M)$ is a $\fg^{x}$-module and $\fg_x$-module.
This gives the functor
$\DS_x: M\mapsto \DS_x(M)$  from the category of $\fg$-modules with semisimple action of $c$ to
the category of $\fg_x$-modules.

There are canonical isomorphisms $\DS_x(\Pi(N))\cong\Pi(\DS_x(N))$ and
$$\DS_x(M)\otimes\DS_x(N)\cong\DS_x(M\otimes N).$$

\subsection{The map $\ds_x$}\label{section-ds-on-red-groth}
Let $\fa$ be any subalgebra of $\fg$.  Write $\fa^x$ for the kernel of $\operatorname{ad}x$ on $\fa$, and similarly for $\fa^c$. We view
$$\fa_x:=\fa^x/([\fg^c,x]\cap \fa^c)$$
as a subalgebra of $\fg_x$.

Let $\CC(\fg)$ be a full subcategory of the category of $\fg$-modules with semisimple action of $c$ and
$\CC(\fa)$ (resp., $\CC(\fg^x),\CC(\fa^x)$) be a full subcategory  of  the category of
$\fa$-modules (resp., $\fg^x,\fa^x$) such that the restriction functors
$$\Res^{\fg}_{\fa}:\CC(\fg)\to\CC(\fa),\ \ \ \Res^{\fg}_{\fg^x}: \CC(\fg)\to \CC(\fg^x),\ \ \ \Res^{\fa}_{\fa^x}:\CC(\fa)\to\CC(\fa^x).$$
are well-defined and that for each $N\in \CC(\fg)$ the $\fg^x$-modules $N^x$ and $xN^c$ 
lie in $\CC(\fg^x)$ (note that $N^x$ and $xN^c$ are submodules of $\Res^{\fg}_{\fg^x}(N)$).

We denote by $\CC(\fg_x)$ (resp., by $\CC(\fa_x)$)
the full category of $\fg_x$-modules $N$ satisfying 
 $\Res^{\fg_x}_{\fg^x}(N)\in \CC(\fg^x)$ (resp.,    $\Res^{\fa_x}_{\fa^x}(N)\in \CC(\fa^x)$).

For $\fm=\fg,\fa,\fg^x,\fa^x,\fg_x,\fa_x$ we denote the reduced Grothendieck group $\Gr_-(\CC(\fm))$ by $R({\fm})$, for ease of notation.

\subsubsection{}
Take $M\in\CC(\fg)$ and set $N:=\Res^{\fg}_{\fg^x}M$. 
The action of $x$ gives a
$\fg^x$-homomorphism $\theta_N: N^c\to \Pi(N^c)$. One has 
$$\theta_N \theta_{\Pi(N)}=0,\ \ \ Im\, \theta_{\Pi(N)}=\Pi( Im\, \theta_N)$$
and
$\DS_x(M)=\Ker\theta_N/Im\,\theta_{\Pi(N)}$ as $\fg^x$-modules.
Using the exact sequences
$$0\to Im\,\theta_{\Pi(N)} \to \Ker \theta_N \to \DS_x(M)
\to 0,\ \ \ \ \
0\to \Ker \theta_N\to N^c \to Im\, \theta_N\to 0
$$
we obtain  
$[\Res^{\fg}_{\fg^x}M^c]=[\DS_x(M)]$ in
$R(\fg^x)$.   However if $M_r$ is the $r\neq0$ eigenspace of $c$ on $M$, $x:M_r\to M_r$ will define an $\fg^x$-equivariant isomorphism of $M_r$, and thus we have $[\Res^{\fg}_{\fg^x}M_r]=0$.  It follows that $[\Res^{\fg}_{\fg^x}M]=[\Res^{\fg}_{\fg^x}M^c]=[\DS_x(M)]$.
Since $\DS_x(M)$ is a $\fg_x$-module this gives the commutative diagram
\begin{equation}\label{eqn-dsx-diagram}
\xymatrix{ R(\fg)\ar[rd]_{\ds_x} \ar[r] & R(\fg^x)  \\
      &        R(\fg_x) \ar[u]}
\end{equation}
where $\ds_x: R(\fg)\to R(\fg_x)$ is given by $\ds_x([M]):=[\DS_x(M)]$,
and the two other arrows are induced by the restriction functors 
$\Res^{\fg}_{\fg^x}$, $\Res^{\fg_x}_{\fg^x}$ respectively.

\subsubsection{Remark}
If $\CC(\fg),\CC(\fg_x)$ are closed under $\otimes$, then
$\ds_x$ is a ring homomorphism.

\subsubsection{Example}
Suppose that $\fg$ is quasireductive.  Recall that $\cF(\fg)$ a  rigid tensor category with the duality $N\mapsto N^*$
given by the anti-automorphism $-\Id$; 
since $\DS_x$ is a tensor functor, it preserves the $*$-duality, so
$$\ds_x: \Gr_-(\fg)\ \to \ \Gr_-(\fg_x)$$
is a ring homomorphism compatible with $*$.  

\subsection{$\ds_x$ and restriction} 

We present results which explain the relationship between $\ds_x$ and the restriction functor.
\begin{lemma}\label{lemma-ds-restriction}
	Suppose that we have a splitting $\fg_x\subseteq\fg^x$ so that $\fg^x=\fg_x\ltimes [x,\fg^c]$.  Then for $M$ in $\mathcal{C}(\fg)$ we have
	\[
	\ds_x[M]=[\Res_{\fg_x}^{\fg}M].
	\]
\end{lemma}
\begin{proof}
	This follows immediately by applying the restriction $R(\fg^x)\to R(\fg_x)$ to our equality $[DS_xM]=[\Res_{\fg^x}^{\fg}M]$.
\end{proof}
	
\begin{lemma}\label{lemma-composition-ds}
	Let $y\in\fg_{\ol{1}}$ with $[y,y]=d$ where $\ad d$ acts semisimply on $\fg$ and $d,c+d$ act semisimply on all modules in $\CC(\fg)$.  Suppose further that $[x,y]=0$, and that we have splittings
	\[
	\fg^{y}\subseteq\fg_y\ltimes [y,\fg^d],\ \ \ \fg^{x+y}=\fg_{x+y}\ltimes[x+y,\fg^{c+d}],
	\]
	Furthermore suppose that under these splittings, $x\in\fg_{y}$ and 
	\[
	(\fg_{y})^{x}=\fg_{x+y}\ltimes[x,\fg_{y}^c].
	\]
	Then we have
	\[
	ds_{x+y}=ds_{x}\circ ds_{y}:R(\fg)\to R(\fg_{x+y})
	\]
\end{lemma}

\begin{proof}
	This follows immediately from Lemma \ref{lemma-ds-restriction} and the corresponding statement for restriction.
\end{proof}

\begin{proposition}\label{proposition-ds-res-diagram}
We have 
the following commutative diagram
$$ 
\xymatrix{ R(\fg) \ar[d]_{\ds_x} \ar[r] & R(\fa^x)  \\
              R(\fg_x)  \ar[r] & R(\fa_x)\ar[u]_{\res^{\fa_x}_{\fa^x}}  }
$$
where the horizontal arrows are induced by  the corresponding restriction functors
and $\res^{\fa_x}_{\fa^x}$ is induced by the morphism $\fa^{x}\to\fa_{x}$.
\end{proposition}
\begin{proof}
The restriction functors give the commutative diagram
$$\xymatrix{ \CC(\fg^x)  \ar[r] & \CC(\fa^x) \\
              \CC(\fg_x)  \ar[r]  \ar[u] & \CC(\fa_x)\ar[u] \\
  }$$
which, in combination with (\ref{eqn-dsx-diagram}) gives the diagram
$$ 
\xymatrix{ R(\fg) \ar[dr]_{\ds_x} \ar[r] & R(\fg^x)\ar[r] & R(\fa^x)  \\
              & R(\fg_x) \ar[u]  \ar[r] & R(\fa_x)\ar[u]  }
$$
where all  arrows except $\ds_x$ are induced by the  restriction functors. By~(\ref{eqn-dsx-diagram}), the above diagram is commutative, and we obtain our result.
\end{proof}

\subsubsection{Example: $\cF(\fg)$ for $\fg$ quasireductive}
Let $E$ be a simple $\fa^x$-module. By Proposition \ref{proposition-ds-res-diagram} a finite-dimensional $\fg$-module $N$ we have

$$\smult(\DS_x(N); E)=\smult(\Res^{\fg}_{\fa^x} N;E).$$

For example, let $\fg$ be a classical Lie superalgebra in the sense of~\cite{KLie} and  $\fa:=\ft$
be a Cartan subalgebra of $\fg_{\ol{0}}$.  The map
$R(\fg)\to R(\ft)$ is given $[N]\mapsto \sch N$.
If $\ft_x$ is  a Cartan subalgebra of $(\fg_x)_{\ol{0}}$, then the composed map
$R(\fg)\to R(\ft_x)$ is given by
$[N]\mapsto \sch \DS_x(N)$. If we fix an embedding $\ft_x\to \ft^x$, we obtain the Hoyt-Reif formula~\cite{HR}
$$\sch \DS_x(N)=(\sch N)|_{\ft_x}.$$

\subsection{A special case}\label{section-ds-a=h}
Consider the case when
$\fg,\fa=\fh$ are as  in Section \ref{section-maxl-qtoral-subalg}.
Denote by $F^x$ the restriction of the form $F$ to $\fh^x$
and set
$$I_x:=\{\lambda\in \ft^*|\ \lambda([x,\fg^c]\cap \ft)=0\}.$$
By assumption, given $\lambda\in I_x$ we have that $\lambda|_{\ft^x}$ lies in the subspace $(\ft_x)^*$.  For $\lambda\in I_x$ we denote this element by $\lambda_x\in(\ft_x)^*$.

Note that for $\nu\in I_x$ satisfying $\rk F^x_{\nu|_{\ft^x}}=\rk F_{\nu}$
the module $\Res^{\fh}_{\fh_x} C(\nu)$ is  simple, so
 $\Res^{\fh}_{\fh_x} C(\nu)\cong\Pi^i C(\nu_x)$ for some $i$.

\begin{corollary}\label{corollary-sch-res-ds_x}
Take $N\in\cF(\fg)$. If
$\sch_{\fh}(N)=\sum\limits_{\nu} m_{\mu}[C(\mu)]$, then
$$\sch_{\fh_x}\bigl(\DS_x(N))\bigr)
=\sum_{\mu\in I_x:\rk F^x_{\mu}=\rk F_{\mu} } m_{\mu} (-1)^{i_{\mu}}\ [C(\mu_x)].$$
where $\Res^{\fh}_{\fh_x} C(\nu)\cong\Pi^{i_{\mu}} C(\mu_x)$.
\end{corollary}
\begin{proof}
Recall that $\sch_{\fh}$ gives an embedding
of $\Gr_-(\fg)$ to $\Gr_-(\fh)$.
Applying Proposition \ref{proposition-ds-res-diagram} to $\cF(\fg)$ we obtain for 
 $R(\fm):=\Gr_-(\fm)$ the commutative diagram 
$$ 
\xymatrix{ R(\fg) \ar[d]_{\ds_x} \ar[r]^{\sch_{\fh^x}} & R(\fh^x)  \\
              R(\fg_x)  \ar[r]^{\sch_{\fh_x}}  & R(\fh_x)\ar[u]_{\res_{\fh^x}^{\fh_x}}  }
$$
where 
 $\res^{\fh^x}_{\fh_x}: R_{\fh_x}\to R^{\fh_x}$ is induced by the map $\fh^x\to\fh_x$. 
In light of Proposition \ref{proposition-cliff-res} (v) we have
$\sch_{\fh^x}(C(\mu))=0$ except for the case
$\rk F^x_{\mu}=\rk F_{\mu}$.
\end{proof}

\subsubsection{Example}
If, in addition,
$(\fg_x)_{\ol{0}}^{\fa_x}=\fa_x$, then, by Corollary \ref{corollary-res-Gr-embedding},
$\sch_{\fa_x}$ gives an embedding
of the reduced Grothendieck ring of $\fg_x$ to $\Gr_-(\fh	_x)$.

For $\fg=\fgl(m|n),\osp(m|n), \fp_n,\fq_n,\fsq_n$  and the exceptional
Lie superalgebras, for each $x$ we can choose a suitable $\fh$
such that $(\fg_x)_{\ol{0}}^{\ft_x}=\ft_x$, see~\cite{DS},\cite{S},\cite{Gcore}.  

\section{The reduced Grothendieck ring for $\cF(\fq_n)$}\label{section-red-Gr-q}
In this section we describe  the  reduced Grothendieck ring for $\cF(\fq_n)$, the category of finite-dimensional $\fq_n$-modules with semisimple action of $(\fq_n)_{\ol{0}}$.  In addition we will explicitly describe the homomorphisms $\ds_s: \Gr_-(\fq_n)\to \Gr_-(\fq_{n-2s})$ induced by the DS functor, to be defined.  Wherever it is not stated, we set $\fg:=\fq_n$.  We will mainly concentrate on the category $\cF(\fq_n)_{int}$ which is the full subcategory of finite-dimensional $\fq_n$-modules with integral weights, and then reduce the corresponding results for $\cF(\fq_n)$  to $\cF(\fq_n)_{int}$. 

\subsection{Structure of $\fq_n$}
Recall that
 $\fq_n$ is the subalgebra of $\fgl(n|n)$ consisting of the matrices with the block form
$$T_{A,B}:=\begin{pmatrix}
A & B\\
B & A
\end{pmatrix}$$

\subsubsection{}
One has $\fg_{\ol{0}}=\fgl_n$. The group $GL_n$ acts on $\fg$ by the inner
automorphisms;
all triangular decompositions of $\fq_n$ are $GL_n$-conjugated.
We denote by $\ft$ the  Cartan subalgebra of $\fgl_n$
spanned by the elements $h_i=T_{E_{ii},0}$ for $i=1,\ldots,n$, where $E_{ij}$ denotes the $(i,j)$ elementary matrix.
 Let  $\{\vareps_i\}_{i=1}^n\subset\ft^*$ 
 be the   basis dual to $\{h_i\}_{i=1}^n$.  The algebra $\fh:=\fq_n^{\ft}$ is
a Cartan subalgebra of $\fq_n$; one has $\fh_{\ol{0}}=\ft$.
The elements $H_i:=T_{0,E_{ii}}$ form a basis of $\fh_{\ol{1}}$; one has
$[H_i,H_j]=2\delta_{ij}h_i$.

 \subsubsection{}
  We  write $\lambda\in\ft^*$
as $\lambda=\sum_{i=1}^n\lambda_i\vareps_i$ and denote by $\Nonzero(\lambda)$ 
the set of non-zero  elements in the multiset $\{\lambda_i\}_{i=1}^n$
and by $\zero\lambda$ the number of zeros in the multiset $\{\lambda_i\}_{i=1}^n$.
Recall that  $\rank F_{\lambda}$
is equal to the cardinality of $\Nonzero(\lambda)$ ($=n-\zero\lambda$).

We call a weight $\lambda\in\ft^*$ {\em integral } (resp., {\em half-integral}) 
if $\lambda_i\in\mathbb{Z}$ (resp., $\lambda_i-\frac{1}{2}\in\mathbb{Z}$) for all $i$. We call a weight $\lambda$  {\em typical} if $\lambda_i+\lambda_j\not=0$ for all $i,j$ and {\em atypical} otherwise; in particular if $\lambda$ is typical then $\zero(\lambda)=0$.

We fix the usual triangular decomposition:
$\fg=\fn^-\oplus\fh\oplus\fn$,
where $\Delta^+=\{\vareps_i-\vareps_j\}_{1\leq i<j\leq n}$.

\subsection{The monoid $\Xi$ and cores}\label{section-monoid-Xi-cores}
We denote by $\Xi$ the set of finite multisets $\{a_i\}_{i=1}^s$
with  $a_i\in\mathbb{C}\setminus\{0\}$ and $a_i+a_j\not=0$ for all $1\leq i,j\leq s$.

We assign to  each finite multiset $A:=\{a_i\}_{i=1}^s$ with  
$a_i\in\mathbb{C}$ the multiset $\Core(A)\in \Xi$,
obtained by throwing out all zeros and the maximal number of elements $a_i,a_j$ with 
$i\not=j$ and $a_i+a_j=0$; for example, 
\[
\Core(\{1,1,-1,-1\}=\emptyset, \ \text{ and } \ \Core(\{1,1,0,0,0,-1\})=\{1\}.
\]

We view $\Xi$ as a commutative monoid with respect to the operation
$$A\diamond B:=\Core(A\cup B)$$
($\emptyset$ is the identity element in $\Xi$).

For $\lambda\in\ft^*$ we set 
\[
\Core(\lambda):=\Core(\{\lambda_i\}_{i=1}^n)=\Core(\Nonzero(\lambda)),
\]
and denote by $\chi_{\lambda}$ the central character of $L(\lambda)$. From~\cite{Serq} it follows that $\chi_{\lambda}=\chi_{\nu}$
 if and only $\Core(\lambda)=\Core(\nu)$.  Observe that $\lambda\in I_0$ if and only if the cardinality of
 $\Core(\lambda)$ is even.

\subsubsection{Relation $\sim$}\label{section-relation-simq}

Recall that we write $\lambda\sim\nu$ if $L(\lambda)$, $L(\nu)$ lie in the same block in the BGG category $\CO$. By above,
$$\lambda\sim\nu\ \ \Longrightarrow\ \ \Core(\lambda)=
\Core(\nu).$$
It is known that the above implication becomes an equivalence if both $\lambda,\nu$ are integral 
or half-integral; in fact it follows from the following fact:
\begin{equation}\label{eqn-simq-from-atyp}
\lambda_i+\lambda_j=0\ \ \Longrightarrow\ \ \lambda-\vareps_i+\vareps_j\sim \lambda
\end{equation}
(this easily follows from the formula for Shapovalov determinants
established in~\cite{Gq}, Thm. 11.1: 
from this formula it follows that for $i<j$ the module
$M(\lambda)$ has a primitive vector of weight $\lambda-\vareps_i+\vareps_j$
if  $\lambda$ is a ``generic weight'' satisfying $\lambda_i+\lambda_j=0$;
the usual density arguments (see~\cite{BGG},\cite{KK}) imply
that $M(\lambda)$ has a primitive vector of weight $\lambda-\vareps_i+\vareps_j$
if  $\lambda$ is any weight satisfying $\lambda_i+\lambda_j=0$).

\subsubsection{Dominant weights}
Recall that $P^+(\fg)$ denotes the set of dominant weights, i.e.
$$P^+(\fg):=\{\lambda\in\ft^*|\ \dim L(\lambda)<\infty\}.$$

By~\cite{Pe}, $\lambda\in P^+(\fg)$ if and only if
 $\lambda_i-\lambda_{i+1}\in\mathbb{N}$ and $\lambda_i=\lambda_{i+1}$
implies $\lambda_i=0$. This implies the following properties:

\begin{enumerate}
\item $P^+(\fg)\cap I_0=P^+(\fg_{\ol{0}})'\cap I_0$; in particular
\[
P^+(\fg)=\{\lambda\in P^+(\fg_{\ol{0}}):\det(w|_{\fh_{\ol{1}}/K_{\lambda}})=1\text{ for all }w\in\Stab_W\lambda\}.
\]
\item if $\lambda\in P^+(\fg)$ is atypical, then $\lambda$ is either integral or half-integral;

\item the set $\Nonzero(\lambda)$ uniquely
determines a dominant weight for a fixed $n$.
\end{enumerate}

\subsubsection{Grading on $C(\nu)$} For $\nu\in\ft^*$, we set
$$T_{\nu}:=H_{i_1}\ldots H_{i_k}$$
where $\Nonzero(\nu)=\{\nu_{i_1}\geq \nu_{i_2}\geq \ldots\geq \nu_{i_k}\}$, and if $\nu_{i_j}=\nu_{i_{j+1}}$ then we require that $i_{j}\geq i_{j+1}$.
This formula determines $T_{\nu}$ uniquely;
for $\nu\in P(\fg)'$ note that we have
$T_{w\nu}=w T_{\nu}$ for each $w\in W$.

Since  $T_{\nu}^2=(-1)^{\frac{k(k-1)}{2}}h_{i_1}\ldots h_{i_k}$  the function 
$t:\ \ft^*\to \{c\in\mathbb{C}|\ c>0\}$
is given by
$$ t(\lambda)^2=(-1)^{\frac{k(k-1)}{2}}\prod_{i: \lambda_i\not=0} \lambda_i,\ \ \text{ where }k:=\rank F_{\lambda};$$
for $\lambda\in I_0$ we obtain 
$ t(\lambda)^2=(-1)^{\frac{\rank F_{\lambda}}{2}}
\prod_{i: \lambda_i\not=0} \lambda_i$.

As explained in Remark \ref{remark-rigidify-parity}, this function $t$ together with a highest weight $\lambda$ determines uniquely irreducible modules $L(\lambda)$ for each $\lambda\in P^+(\fq_n)$.

\subsubsection{Remark} If $\nu\sim0$, then $t(\nu)\in\mathbb{R}^+$.  It follows that if $\lambda\in\ft^*$, $\nu\sim 0$, and $[C(\lambda)][C(\nu)]=\pm[C(\lambda+\nu)]$, then $t(\lambda)t(\nu)=t(\lambda+\nu)$.

\subsubsection{}
Note that $\lambda=-w_0\lambda$ if and only if 
$\Core(\lambda)=\emptyset$; in this case Corollary \ref{corollary-longest-elt-det}
gives $L(\lambda)^*\cong \Pi^{\frac{\rank F_{\lambda}}{2}} L(\lambda)$.

\subsubsection{Example}\label{example-multiples-alpha}
For $\alpha=\vareps_i-\vareps_j$ one has $T_{\alpha}=H_jH_i$. Since for
$e\in \fg_{\alpha}\cap \fg_{\ol{0}}$ 
one has $(\ad H_j)(\ad H_i) e=e$, we obtain
$\fg_{\alpha}\cong C_\alpha$. Moreover, 
$S^k C(\alpha)= C(k\alpha)$  for $\alpha\in\Delta$.

\subsection{Embedding into exterior algebra}
Let $\fh^{\mathbb{Z}}$ be the Lie subalgebra of $\fh$ over $\mathbb{Z}$
generated by $H_1,\ldots,H_n$, and let $\cU(\fh^{\mathbb{Z}})$ be the integral enveloping algebra of $\fh^{\mathbb{Z}}$.
Consider the canonical epimorphism $\phi_0: \cU(\fh^{\mathbb{Z}})\to \cS(\fh^{\mathbb{Z}}_{\ol{1}})$, where
$\cS(\fh^{\mathbb{Z}}_{\ol{1}})$ is the exterior ring generated by $\xi_i:=\phi_0(H_i)$.
Note that $\cS(\fh^{\mathbb{Z}}_{\ol{1}})$ is a $\mathbb{Z}$-graded supercommutative ring, free over $\mathbb{Z}$ with basis $\xi_{i_1}\cdots\xi_{i_j}$ with
$1\leq i_1<i_2<\ldots<i_j\leq n$. 

Let $\CC_{\mathbb{R}}$ be the subcategory of $\cF(\fh)$ consisting of weights $\lambda$ such that $\lambda_i\in\mathbb{R}$ for all $i$.  In particular, $\frac{t_{\lambda}}{|t_{\lambda}|}\in\{1,\sqrt{-1}\}$.

\subsubsection{}
We view the ring
$B:=\mathbb{Z}[e^{\nu}:\nu\in\ft^*]\otimes_{\mathbb{Z}} \cS(\fh_{\ol{1}}^{\mathbb{Z}})$ as a $\mathbb{Z}$-graded supercommutative ring by defining the degree of  $\mathbb{Z}[e^{\nu}:\nu\in\ft^*]$ to be zero.  We construct the ring $B^{spoil}$ as in~\ref{section-spoiled-algs}.

\begin{proposition}\label{proposition-red-Gr-embedding-spoiled-polys}
	The map $[C(\lambda)]\mapsto \vareps^i\frac{t(\lambda)}{|t(\lambda)|}\cdot e^{\lambda} \phi_0(T_{\lambda})\ \text{ for }\lambda\in I_i$
gives a ring monomorphism
\[
\Gr_-(\CC_{\mathbb{R}})\hookrightarrow B^{spoil}\otimes_{\mathbb{Z}} \mathbb{Z}[\sqrt{-1}]
\]
 which is compatible with  the action of $W$. 
One has
	$$\phi_0(T_{\lambda})=\xi_{i_1}\ldots\xi_{i_k},$$
	where $\Nonzero(\lambda)=\{\lambda_{i_1}\leq \ldots\leq \lambda_{i_k}\}$.
\end{proposition}

\begin{proof} 
	 That the map is injective and $W$-equivariant is straightforward.  We show that it is an algebra homomorphism.  First observe that:
\[
[C(\nu)][C(\lambda)]\not=0\ \Longleftrightarrow\   [C(\nu)][C(\lambda)]= \pm [C(\lambda+\nu)]\ \Longleftrightarrow  \ \phi_0(T_{\lambda}T_{\nu})\not=0
\ \Longleftrightarrow  \ T_{\lambda}T_{\nu}=\pm T_{\lambda+\nu}.
\]
Suppose that $[C(\nu)][C(\lambda)]\not=0$. Then we may write
$T_{\lambda}T_{\nu}=(-1)^j T_{\lambda+\nu}$ and $\phi_0(T_{\lambda})\phi_0(T_{\nu})=(-1)^{j}\phi_0(T_{\lambda+\nu})$.

Fix even vectors $v_{\lambda}\in C(\lambda)_{\ol{0}}$ and $v_{\nu}\in C(\nu)_{\ol{0}}$.
Since $H_i C(\lambda)=0$ if $\lambda_i=0$ we have
\[
T_{\lambda} T_{\nu}(v_{\lambda}\otimes v_{\nu})=T_{\lambda}v_{\lambda}\otimes T_{\nu}v_{\nu}=\frac{t(\lambda)t(\nu)}{|t(\lambda)t(\nu)|}v_{\lambda}\otimes v_{\nu}.
\]
Note that $v_{\lambda}\otimes v_{\nu}$ is an even vector of $\Pi^i C(\lambda+\nu)\cong C(\lambda)\otimes C(\nu)$.  On the other hand we have
\[
T_{\lambda} T_{\nu}(v_{\lambda}\otimes v_{\nu})=(-1)^{j} T_{\lambda+\nu}(v_{\lambda}\otimes v_{\nu})=(-1)^{j+i} t(\lambda+\nu) v_{\lambda}\otimes v_{\nu}.
\]
Thus we have that $\frac{t(\lambda+\nu)}{|t(\lambda+\nu)|}=(-1)^{i+j}\frac{t(\lambda)t(\nu)}{|t(\lambda+\nu)|}$.  From the above equalities it is easy to check our map is a homomorphism.  Injectivity is straightforward.
\end{proof}

\subsubsection{Remark} If we extend scalars to $\mathbb{C}$, the above map defines an embedding
\[
\Gr(\fh)_{\mathbb{C}}\hookrightarrow B_{\ol{0}}\otimes_{\ZZ}\mathbb{C}.
\]
However this kills the two-torsion part of $\Gr(\fh)$.

\subsubsection{}\label{subsection-core-diamond-reln}
One has 
$$\phi_0(T_{\lambda}T_{\nu})\not=0\ \Longleftrightarrow \ T_{\lambda}T_{\nu}= \pm T_{\lambda+\nu}\ \Longrightarrow\  \Core(\lambda+\nu)=\Core(\lambda)\diamond \Core(\nu).$$
For example, for $[C(2,0,0,1)][C(0,-2,-1,0)]=\pm [C(2,-2,-1,1)]$ and 
$$\Core(\lambda)=\{2,1\},\ \ \Core(\nu)=\{-2,-1\},\ \  
\Core(\lambda+\nu)=\emptyset.$$

\subsubsection{}
Recall that $\Gr_-(\fh)$ has a finite $\mathbb{Z}$-grading (see
\ref{eqn-grading-red-Gr-h}). By Section \ref{subsection-core-diamond-reln},  $\Gr_-(\fh)$ is also $\Xi$-graded: set $\Gr_-(\fh)_{A}$ to be spanned by $[C(\lambda)]$ with $\Core(\lambda)=A$.  Then
\[
\Gr_-(\fh)=\bigoplus_{A\in \Xi} \Gr_-(\fh)_{A}\ \text{ with } \ \Gr_-(\fh)_{A} \Gr_-(\fh)_{B}\subset \Gr_-(\fh)_{A\diamond B},
\]
Note that $\Gr_-(\fh)_A=0$ if the cardinality of $A$ is odd or greater than $n$.  Further $\Gr_-(\fh)_A$ is $W$-stable.

\begin{corollary}\label{corollary-red-Gr-ring-struc}
\begin{enumerate}
\item The subring $\Gr_-(\fh)_{\emptyset}\cap \Gr_-(\fh)_{int}$ is spanned by
$[C(\nu)]$ with $\nu\sim 0$. This subring is generated by $[C(k\alpha)]\ $
for $\alpha\in\Delta$ and $k\in\mathbb{Z}$.

\item
If $\nu\sim 0$ and
$[C(\nu)][C(\lambda)]\not=0$, then $[C(\nu)][C(\lambda)]=\pm  [C(\lambda+\nu)]$ and $\lambda+\nu\sim \lambda$.
\end{enumerate}
\end{corollary}

\begin{proof}
The assertion (i) follows from Section \ref{section-relation-simq}.
By (i) it is enough to verify (ii) for $\nu:=k(\vareps_i-\vareps_j)$. In this case
the inequality $[C(\nu)][C(\lambda)]\not=0$  implies $\lambda_i=\lambda_j=0$.
By~(\ref{eqn-simq-from-atyp}) for such $\lambda$ one has $\lambda+k(\vareps_i-\vareps_j)\sim \lambda$.
\end{proof}

\subsubsection{Remark} All integral, non-typical blocks $\cB$ of $\fq_n$ have that $\Pi\cB=\cB$ (indeed, if $\cB$ is integral and not typical then it admits a simple module $L(\lambda)$ such that $\zero(\lambda)>0$; now we conclude by Thm.~4.1 of \cite{Nicki}).  Thus we may consider for such blocks the corresponding reduced Grothendieck group $\Gr_-(\cB)$.  Then the embedding in Proposition \ref{proposition-red-Gr-embedding-spoiled-polys} exists over the coefficient ring $\mathbb{Z}[\sqrt{-1}]$.   If we let $\cB_0$ denote the principal block of $\fq_n$, then $\Gr_-(\cB_0)=\Gr_-(\fh)_{\emptyset}\cap \Gr_-(\fh)_{int}$ is a ring.   In this case, because $t(\lambda)\in\mathbb{R}^+$ for all $\lambda\sim0$, the map in Proposition \ref{proposition-red-Gr-embedding-spoiled-polys} descends to an embedding with {\em integral} coefficients for $\Gr_-(\cB_0)$.

\subsection{Supercharacters of some highest weight modules}

\begin{lemma}\label{lemma-sch-Un-}
$\sch_{\fh} M(\lambda)=\displaystyle\sum_{\nu\sim \lambda} m_{\nu} [C(\nu)]$.
\end{lemma}
\begin{proof}
We start from $\lambda=0$; in this case
$\sch_{\fh} M(0)=\sch_{\fh}\mathcal{U}(\fn^-)$.
Let $\alpha_1,\dots,\alpha_N$ denote the negative roots of $\fg$.  We have
	\begin{equation}\label{eqn-Un-wt-space}
	\mathcal{U}(\fn^-)_{\nu}=\bigoplus\limits_{\sum k_j\alpha_j=\nu}\bigotimes_j S^{k_j}\fg_{\alpha_j}.
	\end{equation}
By Example \ref{example-multiples-alpha},
 $S^{k}\fg_{\alpha}\cong C(k\alpha)$.
Since $k\alpha\sim 0$,
the required formula follows Corollary \ref{corollary-red-Gr-ring-struc} (i).

For an arbitrary $\lambda$ one has
$\sch_{\fh} M(\lambda)=\sch_{\fh}\mathcal{U}(\fn^-)\otimes [C(\lambda)]$
and the result follows from Corollary \ref{corollary-red-Gr-ring-struc} (ii).
\end{proof}

We need some terminology for the following proposition: for roots $\alpha_1=\epsilon_{i_1}-\epsilon_{j_1},\alpha_2=\epsilon_{i_2}-\epsilon_{j_2}$, we write $\alpha_1\prec\alpha_2$ if $\max(i_1,j_1)<\min(i_2,j_2)$.  If $\alpha=\epsilon_i-\epsilon_j$ is any root, we call the set $\{i,j\}$ its support. We note that if we have roots $\alpha_{1},\dots,\alpha_j$ with non-overlapping supports, and positive integers $k_1,\dots,k_j$, then $S^{j_1}\fg_{\alpha_1}\otimes\cdots\otimes S^{k_j}\fg_{\alpha_j}$ is an irreducible $\fh$-module. 

\begin{proposition}
	We have
	\[
	\sch_{\fh}\mathcal{U}(\fn^-)=\sum\limits_{\nu\sim 0}m_{\nu}[C(\nu)]
	\]
	where each coefficient $m_{\nu}$ is either $0$ or $\pm1$.
	
	Further if $m_{\nu}\neq0$, then $\nu$ may uniquely be written as
	\[
	\nu=k_1(\alpha_{11}+\dots+\alpha_{1j_1})+\dots+k_r(\alpha_{r1}+\dots+\alpha_{rj_r}),
	\]
	where all of $k_1,\dots,k_r$ are distinct, the supports of all $\alpha_{ij}$ are distinct, and $\alpha_{i1}\prec\dots\prec\alpha_{ij_{i}}$.  In this case we have
	\[
	m_{\nu}[C(\nu)]=\prod_{i=1}^{r}[S^{k_i}\fg_{\alpha_{i1}}]\cdots[S^{k_i}\fg_{\alpha_{ij_i}}].
	\]
\end{proposition}

\begin{proof}	
We use the embedding of Proposition \ref{proposition-red-Gr-embedding-spoiled-polys}.  Let $y_{\nu}$ be the image of $[(\cU(\fn^-))_{\nu}]$ in $R(\ft)\otimes \cS(\fh_{\ol{1}})_{\ol{0}}$ (see Proposition \ref{proposition-red-Gr-embedding-spoiled-polys}). Since $t_{k\alpha}=k$ for $\alpha\in \Delta$, by~(\ref{eqn-Un-wt-space}) one has
$$y_{\nu}=e^{\nu} \displaystyle\sum_{(k_1,\ldots,k_N): \sum k_i\alpha_i=\nu} 
\phi_0(T_{k\alpha_i})$$ 
Recall that $\phi_0(T_{\vareps_p-\vareps_q})=\xi_q\xi_p$. In particular,
$\phi_0(T_{\alpha}T_{\beta})=0$ if $(\alpha|\beta)\not=0$
(where $(-|-)$ stands for the usual form on $\ft^*$).
Hence
$$\begin{array}{l}
y_{\nu}=e^{\nu} \displaystyle\sum_{(k_1,\ldots,k_N)\in U } 
 \phi_0(T_{k\alpha_i}),\\ 
\text{ where }
U:=\{(k_1,\ldots,k_N)|\ \sum k_i\alpha_i=\nu,\ 
\ k_ik_j=0 \text{ for }(\alpha_i|\alpha_j)\not=0\}.\end{array}$$

If for any $(k_1,\ldots,k_N)\in U$ we have that $k_i\neq k_j$ for all nonzero $k_i,k_j$ with $i\neq j$ then it is clear $U$ is a singleton set and we are done.

Thus suppose that $(k_1,\dots,k_n)$ has $k_i=k_j\not=0$ for some $i\neq j$, and WLOG suppose $i=1$, $j=2$. Write $\alpha_i=\vareps_{p_i}-\vareps_{q_i}$ for $i=1,2$; 
then $p_i>q_i$ and, since $(\alpha_1|\alpha_2)=0$,
the numbers
$p_1,q_1,p_2,q_2$ are pairwise distinct. If $p_1>q_2$ and $p_2>q_1$
then we have negative roots 
$$\alpha'_1:=\vareps_{p_1}-\vareps_{q_2},\ \ \alpha'_2:=\vareps_{p_2}-\vareps_{q_1}$$
with $\alpha_1+\alpha_2=\alpha'_1+\alpha'_2$. We may assume that $\alpha'_1=\alpha_3$ and $\alpha'_2=\alpha_4$. Since
$(\alpha_j|\alpha_1), (\alpha_s|\alpha_1)\not=0$, 
one has $k_3=k_4=0$, so 
$$( k_1,\ldots,k_N)=(k_1,k_1,0,0,k_5,\ldots,k_N).$$

Observe that $(0,0,k_1,k_1, k_5,\ldots,k_N)\in U$. One has
$$\phi_0(T_{\alpha_1}T_{\alpha_2}+T_{\alpha'_1}T_{\alpha'_2})=
\xi_{p_1}\xi_{q_1}\xi_{p_2}\xi_{q_2}+\xi_{p_1}\xi_{q_2}\xi_{p_2}\xi_{q_1}=0.$$
Therefore we can substitute $U$ by a smaller set, where $k_i=k_j$ implies that 
$\alpha_i=\vareps_{p_i}-\vareps_{q_i}$, $\alpha_j=\vareps_{p_j}-\vareps_{q_j}$ 
are such that $p_i>q_i>p_j>q_j$ or  $p_j>q_j>p_i>q_i$.  From this the result follows.

\end{proof}

\begin{corollary}
One has
$\sch_{\fh} M(\lambda)=\displaystyle\sum_{\nu\sim \lambda} m_{\nu} [C(\nu)]$ with
$m_{\nu}\in \{0,\pm 1\}$.
\end{corollary}

\begin{theorem}\label{theorem-sch-irrep-q} 
	For $\lambda\in\ft^*$ we have $\sch_{\fh}L(\lambda)=\displaystyle\sum_{\nu\sim \lambda} k_{\nu} [C(\nu)]$.
\end{theorem}

\begin{proof}
	With the help of Lemma \ref{lemma-sch-Un-}, the proof works in the exact same fashion as Theorem \ref{theorem-parity_invce}.
\end{proof}

\begin{corollary}
Let $L(\nu)$ be a finite-dimensional module. Then
\[
\smult (L(\lambda)\otimes L(\nu):L(\mu))\neq0\ \ \Longrightarrow\ \ 
\Core(\mu)=\Core(\lambda)\diamond \Core(\nu)\
\]
Note further that when $\smult (L(\lambda)\otimes L(\nu):L(\mu))\neq0$, then either $\lambda\in I_0\ \text{ or }\nu\in I_0$.
\end{corollary}

\subsubsection{Remark}
The Kac-Kazhdan modification of $\CO$
introduced in~\cite{KK} is closed under tensor product; 
the modules in this category are not always of finite length,
but the multiplicity is well-defined (see~\cite{DGK}). 
In this category the above formula holds for arbitrary $L(\lambda)$ and $L(\nu)$.

\subsubsection{}
Recall that $P(\fg):=P^+(\fg)+\mathbb{Z}\Delta$ and for $\nu\in P^+(\fg)$ we set
$$
a_{\nu}:=\sum_{w\in W/\Stab_W\nu} w[C(\mu)].
$$

\begin{lemma}
	Suppose that $|\Nonzero(\nu)|$ is odd with $\nu_i=\nu_j\neq0$ for some $i\neq j$; then for $\lambda\in\ P^+(\fq_n)$, $[L(\lambda)_{\nu}:C(\nu)]$ is even.
\end{lemma}

\begin{proof}

Assume that this does not hold. WLOG we may assume $i=1,j=2$, and set $\alpha=\epsilon_1-\epsilon_2$.  Let $\fh'$ to be the subalgebra of $\fh$ generated by $H_{3},\dots,H_n$, and set $\fq_2(\alpha)$ to be the natural subalgebra of $\fq_n$ isomorphic to $\fq_2$ with weight $\epsilon_1,\epsilon_2$.  Clearly $\fq_2(\alpha)\times\fh'$ is a subalgebra of $\fq_n$.  

Then by Theorem \ref{theorem-sch-irrep-q},
$\Core(\nu)=\Core(\lambda)$. We view
\[
N:=\sum_{i\in\mathbb{Z}}  L(\lambda)_{\nu+i\alpha}
\]
as a $\fq_2(\alpha)\times \fh'$-module. We will write  $\nu=k(\vareps_1+\vareps_2)+\nu'$ where $\nu'$ is the corresponding $\ft'=\fh'_{\ol{0}}$-weight.  We assume that $k>0$, with the case of $k<0$ being similar.  By the representation theory of $\fq_2$, the only irreducible $\fq_2$-modules with $k(\epsilon_1+\epsilon_2)$ are those of the form $L_{\fq_2(\alpha)}(k+i;k-i)$ for some $i>0$.  These are always typical, and are isomorphic to their parity shifts if and only if $i=k$.  

Therefore since $\rank F_{\nu}$ is odd one has
\[
L_{\fq_2(\alpha)}(k+i;k-i)\boxtimes L_{\fh'}(\nu')=\left\{\begin{array}{lc}
L_{\fq_2(\alpha)\times \fh'}(\nu+i\alpha) &\text{ if }  i\neq k\\
L_{\fq_2(\alpha)\times \fh'}(\nu+k\alpha)\oplus \Pi L_{\fq_2(\alpha)\times \fh'}(\nu+k\alpha) & i=k.
\end{array}\right.
\]
Because $L_{\fq_2(\alpha)}(k+i;k-i)$ is typical,
\[
\dim L_{\fq_2(\alpha)\times \fh'}(\nu+i\alpha)_{\nu}=\left\{\begin{array}{lcl}
2\dim C(\nu) &\text{if } & i\neq k;\\
\dim C(\nu) & \text{if } & i=k.
\end{array}\right.
\]
Therefore in order to prove that $\dim N_{\nu}$ is divisible by $2\dim C(\nu)$
it is enough to show that 
$\mult(N; L_{\fq_2(\alpha)\times \fh'}(\nu+k\alpha))$ is even, where $\mult$ denotes that non-graded multiplicity.
By above, for $i\neq k$ one has
\[
\dim L_{\fq_2(\alpha)\times \fh'}(\nu+i\alpha)_{\nu+k\alpha}=2\dim C(\nu+k\alpha).
\]
Hence it is enough to verify that 
$\dim N_{\nu+k\alpha}$ is divisible by $2\dim C(\nu+k\alpha)$ that is
\[
\mult(L(\lambda);C(\nu+k\alpha))\equiv 0 \mod 2.
\]

Let $\{\nu_i\}_{i=1}^n$ contain
$i_+$ copies of $k$ and $i_-$ copies of $-k$.  Then $\{\nu+k\alpha\}_{i=1}^n$ contains
$j_+=i_+-2$ copies of $k$ and $j_-$ copies of $-k$.
Therefore $i_+-i_-\not=j_+-j_-$, so 
$$\Core(\nu+k\alpha)\not=\Core(\nu)=\Core(\lambda)$$ 
and thus by Theorem \ref{theorem-sch-irrep-q}, $\mult(L(\lambda);C(\nu+k\alpha)\equiv 0 \mod 2$
as required.

\end{proof}

Combining Theorem \ref{theorem-qred-sch}, Theorem \ref{theorem-sch-irrep-q}, and the results of \cite{GS}, we obtain our main result:

\begin{theorem}\label{theorem-sch-q-long}
	$\sch_{\fh}$ defines a morphism of spoiled superalgebras $\Gr_-(\fg)\to\Gr_-(\fh)^W$.  Further:
\begin{enumerate}
	\item For $\lambda\in P^+(\fq_n)\cap I_0$ one has
	\[
	\sch_{\fh} L(\lambda)=\sum_{\substack{\nu\in P^+(\fg)\cap I_0:\ \nu\sim\lambda \\ \zero\lambda\geq \zero\nu, \ 
		\zero\lambda-\zero \nu\equiv 0\ \mod 4}} k_{\nu} a_{\nu},\ \ k_{\nu}\in\mathbb{N}.
	\]
	\item  For $\lambda\in P^+(\fg)\cap I_1$ one has
	\[
	\sch_{\fh} L(\lambda)=\sum_{\substack{\nu\in P^+(\fg)\cap I_1:\ \nu\sim\lambda \\ \zero\lambda\geq \zero\nu}} k_{\nu} a_{\nu},\ \ k_{\nu}\in\{0,1\}.
	\]
\end{enumerate}
In both cases $k_{\lambda}=1$, and $\sch_{\fh} L(\lambda)=a_{\lambda}$  if $\lambda\in P^+(\fg)$ is typical. 
\end{theorem}

\begin{proof}
	The only part that remains to be justified is the inequality $\zero\lambda\geq\zero\nu$.  For this we invoke the results of \cite{GS}, where it is shown that if $x=H_1\in(\fq_{n})_{\ol{1}}$, and $\lambda\in P^+(\fq_n)$, then $ds_x^{\zero\lambda+1}L(\lambda)=0$.  
\end{proof}

\begin{corollary}
	For $\lambda\in P^+(\fg)$ and $\nu\notin WP^+(\fg)$ we have $2\dim E_{\nu}$ divides $\dim L(\lambda)_{\nu}$.
\end{corollary}

\begin{corollary}\label{corollary-red-Gr-q(n)}
\begin{enumerate}
\item
The map $N\mapsto \sch_{\fh} N$ induces an isomorphism
$$\Gr_-(\fq_n)_{\mathbb{Q}}\iso \Gr_-(\CC)_{\mathbb{Q}}^W.$$
\item For $A\in\Xi$, let $\cF(\fq_n)_A$ be the full subcategory of $\cF(\fq_n)$ consisting of modules of central character corresponding to $A$.  Then $\sch_{\fh}$ restricts to an isomorphism $\Gr_-(\cF(\fq_n)_A)_{\mathbb{Q}}\to (\Gr_-(\fh)_{A})^W_{\mathbb{Q}}$.
\item 
The duality is given by $a_{\nu}^* =a_{-w_0(\nu)}$
\end{enumerate}
\end{corollary}
\begin{proof}
It is clear that $\sch_{\fh}$ has image lying in $\Gr_-(\CC)_{\mathbb{Q}}^W$, and by Corollary \ref{corollary-res-Gr-embedding}, $\sch_{\fh}$ is an embedding. 
From Theorem \ref{theorem-sch-q-long} and the fact that for each $\lambda\in P^+(\fg)$ the number of elements $\nu\in P^+(\fg)$ satisfying $\nu<\lambda$ and $\nu\sim\lambda$ is finite,
it follows that  the image of $\Gr_-(\fg)$ contains $a_{\nu}$
for each $\nu\in P^+(\fg)$. Hence the image is equal to $\Gr_-(\fh)_{\mathbb{Q}}^W$, proving surjectivity giving (i).  Part (ii) is an easy consequence of (i).

For (iii) we simply apply Theorem \ref{theorem-sch-irrep-q}, and for (iii) we observe that $T_{-\lambda}=\sigma(T_{\lambda})$ and $t(\lambda)=t_{-\lambda}$.  Thus it is easy check that $C(\lambda)^*\cong C(-\lambda)$.
\end{proof}

\subsection{The map $\ds_s$}\label{section-ds-rk-s-q(n)}

For $s=\frac{1}{2},\dots,\frac{n}{2}$, set 
\[
x_s=H_{n+1-2s}+\dots+H_n.
\]
Then $c=x_s^2$ is a semisimple element of $\fg_{\ol{0}}$, and we have $DS_x\fq_n=\fq_{n-2s}$.  Further we have a splitting $\fq_n^c=\fq_{n-2s}\ltimes[x,\fq_n^c]$, where $\fq_{n-2s}\sub\fq_n$ is the natural embedding such that $\ft_x\sub\ft$ is spanned by $h_1,\dots,h_{n-2s}$ and $\Delta(\fg_x)=\{\vareps_i-\vareps_j\}_{1\leq i\not=j\leq n-2s}$.  

We write $DS_s:=DS_{x_s}$ and $ds_s:R(\fq_n)\to R(\fq_{n-2s})$ for the induced homomorphism on reduced Grothendieck rings.  These splittings of $\fg_{x_s}$ in $\fg^{x_s}$ satisfy the hypotheses of Lemma \ref{lemma-composition-ds}, thus we have $ds_{i}\circ ds_{j}=ds_{i+j}$.

\subsubsection{Remark}
	For $i=1,\ldots, \lfloor\frac{n}{2}\rfloor$ set
$$\alpha_i:=\vareps_{n-2s+1}-\vareps_{n-2s+2},\ \ \ \alpha_i^{\vee}:=h_{n-2s+1}-h_{n-2s+2}.$$
Let $y_{\alpha}$ be a nonzero odd element in $\fg_{\alpha}$. For each $1\leq s\leq \frac{n}{2}$, consider
$y_s:=\sum_{i=1}^s x_{\alpha_i}$.  Then clearly $[y_s,y_s]=0$, and one can show that $ds_{y_s}=ds_{2s}$ as defined above.

We view $\ft_{x_s}^*$ as a subspace in $\ft^*$  via the natural embedding 
$\iota_{n,s}: \ft_{x_s}^*\hookrightarrow \ft^*$ (given by $\vareps_i\mapsto \vareps_i\in\ft^*$ for $i=1,\ldots,n-2s$).

\begin{corollary}\label{cor-ds-smult-x_s}
Take $N\in\cF(\fg)$. For each $\nu\in\ft^*_{x_s}$   we have
$$\smult(\DS_s(N):C(\nu))=\smult(N:C(\iota_{n,s}(\nu))).$$
\end{corollary}
\begin{proof}
Retain notation of Section \ref{section-ds-a=h} and set $\mu:=\iota_{n,s}(\nu)\in\ft^*$.
Note that $\ft_{x_s}$ is spanned by $h_1,h_2,\ldots,h_{n-2s}$ and
$\ft^{x_s}=\ft$.

Since $\ft\cap [x_s,\fg^c]$ is spanned by $h_{n-s+1},\dots,h_n$,
for $\mu\in I_x$ 
the restriction of $F_{\mu}$ to $\ft^{x_s}$
written with respect to the above basis has the diagonal entries
$\mu_1,\ldots,\mu_{n-s}$ and zeros on the last $s$ places.
In particular, for $\mu\in I_x$ one has $\rk F_{\mu}^{x_s}=\rk F_{\mu}$ 
if and only if $\mu_{i}=0$ for $i>n-s$, i.e. $\mu\in \ft_x^*$.
By Corollary \ref{corollary-sch-res-ds_x} we obtain
$$\smult(\DS_s(N);C(\nu))=(-1)^{i_{\mu}}\smult(N;C(\mu))$$
where $\Res^{\fh}_{\fh_{x_s}} C(\mu)\cong \Pi^{i_{\mu}}C(\nu))$.
The formulae 
$$T_{\mu}=\prod_{i:\ \mu_i\not=0} H_i=\prod_{i:\ \nu_i\not=0} H_i=T_{\nu},\ \ \ 
t(\mu)=t(\nu)$$
 give $\Res^{\fh}_{\fh_{x_s}} C(\mu)\cong C(\nu)$
as required.
\end{proof}

\begin{corollary}\label{cor-ds-q(n)}
The map $\ds_s: \Gr_-(\fg)\to \Gr_-(\fg_{x_s})$ is given by
$$\ds_s(a_{\mu})=\left\{\begin{array}{lcl}
0 & & \text{ if }
\zero \mu<s\\
a_{\mu'} & & \text{ if }
\zero \mu\geq s\end{array}\right.$$
 where $\mu'\in P^+(\fg_x)$
is such that $\Nonzero(\mu)=\Nonzero(\mu')$. 
\end{corollary}
\begin{proof}
Write $\sch_{\fh} N=\sum_{\mu\in P^+(\fg)} k_{\mu}a_{\mu}$
and $\sch_{\fh_x} \DS_s(N)=\sum_{\nu \in P^+(\fg_x)} m_{\nu}a_{\nu}$.
By Corollary \ref{cor-ds-smult-x_s} $k_{\mu}=m_{\nu}$ if $\mu=\iota_{n,s}(\nu)$.
\end{proof}

\subsubsection{}
We denote by $\cF(\fg)_{int}$ the full subcategory
of $\cF(\fg)$ with the modules whose weights with nonzero weight spaces lie in the lattice generated by $\epsilon_1,\dots,\epsilon_n$. 

\begin{corollary}\label{cor-ds-kernel}
The kernel of
$\ds_s: \Gr_-(\fg)\to \Gr_-(\fg_x)$ is spanned by $a_{\mu}$ with  $\zero \mu<s$ and
the image of $\ds_s$ is equal to $\Gr_-(\cF(\fg_{x_s})_{int})$.
\end{corollary}

\begin{corollary}
For $\lambda\in P^+(\fg)$ and $\nu\in P^+(\fg_{x_s})$ one has
$$\zero\lambda-\zero\nu-s\not\equiv 0 \mod 4\ \ \Longrightarrow\ \ \smult(\DS_s(L(\lambda)):L_{\fg_x}(\nu))=0.$$
\end{corollary}
\begin{proof}
Combining Theorem \ref{theorem-sch-q-long} and Corollary \ref{cor-ds-q(n)} we conclude that
$\ds_s(\sch L(\lambda))$ lies in the span of
$a_{\nu}$ with $\nu\in P^+(\fg_x)$ such that 
$n-zero(\lambda)-(n-s-\zero (\nu))\equiv 0\mod 4$ that is
$\zero\lambda-\zero\nu-s\equiv 0 \mod 4$.

Let $\nu_0$ be maximal (with respect to the standard partial order in $\ft_x^*$)
such that 
$$\smult(\DS_s(L(\lambda));L_{\fg_x}(\nu))\not=0 \ \ \text{ and }
\zero\lambda-\zero\nu-s\not\equiv 0 \mod 4.$$
The maximality of $\nu_0$ forces $\zero\nu\not\equiv \zero\nu_0$ if $\nu>\nu_0$ and
  $L_{\fg_x}(\nu)$ is a subquotient of $\DS_s(L_{\fg}(\lambda))$.
By Theorem \ref{theorem-sch-q-long} we obtain
$[L_{\fg_x}(\nu):C(\nu_0)]=0$
if $\nu\not=\nu_0$ and $L_{\fg_x}(\nu)$ is a subquotient of $\DS_s(L_{\fg}(\lambda))$.  Therefore 
$\smult(\DS_s(L(\lambda));L_{\fg_x}(\nu_0))$ is equal to the coefficient of
$a_{\nu_0}$ in $\sch_{\fh} \DS_s(L(\lambda))=\ds_s(\sch_{\fh} L(\lambda))$,
which is zero by above.
\end{proof}

\subsection{Example: $\fq_2$}
Recall that the atypical dominant weights for $\fq_2$ are of the form
$s(\vareps_1-\vareps_2)$ for $s\in\frac{1}{2}\mathbb{N}$.

\begin{proposition}
Take $\fg=\fq_2$ with $\Delta^+=\{\alpha\}$. For $\lambda\in P^+(\fq_2)$   one has
$$\sch_{\fh} L(\lambda)=\left\{
\begin{array}{lcl}
a_{\lambda} & & \text{ if $\lambda$ is typical};\\
\sum_{i=1}^s a_{i\alpha} & & \text{ if }\lambda=s\alpha,\ s\in\mathbb{N};\\
\sum_{i=0}^s a_{i+\frac{1}{2}\alpha} & & \text{ if }\lambda=(s+\frac{1}{2})\alpha,\ s\in\mathbb{N}.
\end{array}
\right.$$
\end{proposition}

\begin{proof}
If $\lambda$ is typical, the assertion follows from Theorem \ref{theorem-parity_invce}. 
Consider the case when $\mu\in P^+(\fq_2)$ is atypical and $\mu\not=0$.  Write $K(\mu)$ for the maximal finite-dimensional quotient of $\Ind_{\fb}^{\fq_2}C(\mu)$.  It is known that $\sch_{\fh}K(\mu)=a_{\mu}$.  Further, if $\mu=s\alpha$ for $s\in\mathbb{N}$, then it is known (for example see Section 7 of \cite{Gq}) that we have short exact sequences
\[
0\to V_0\to K(\alpha)\to L(\alpha)\to 0, \ \ \ \ 0\to \Pi L((s-1)\alpha)\to K(s\alpha)\to L(s\alpha)\to 0
\]
where $V_0$ is a nontrivial extension of $\mathbb{C}$ and $\Pi\mathbb{C}$.  On the other hand if $\mu=(s+\frac{1}{2})\alpha$ then for $s=0$ we have $K(\alpha/2)=L(\alpha/2)$ and for $s>0$ we have short exact sequences
\[
0\to \Pi L((s-\frac{1}{2})\alpha)\to K((s+\frac{1}{2})\alpha)\to L((s+\frac{1}{2})\alpha)\to 0.
\]
From these results we can obtain the desired formulas by induction.
\end{proof}

\subsection{Realization of $\Gr_-(\fq_n)$}
Recall that $\cF(\fg)_{int}$ denotes the full subcategory of $\cF(\fg)$ consisting
of the modules with integral weights and that $\CC$ denotes the full subcategory of $\cF(\fh)$ consists of modules with weights lying in $P(\fg)$. We write $\CC_{int}$ for the subcategory of $\CC$ consisting of those modules with integral weights.

\subsubsection{} By Theorem \ref{theorem-sch-q-long},
 the map $\sch$ gives embeddings
\[
\Gr_-(\fg)\hookrightarrow \Gr_-(\CC)^W,\ \ \ 
\Gr_-(\cF(\fg)_{int})\hookrightarrow \Gr_-(\CC_{int})^W.
\]
 Further, we have identified the image of $\Gr_-(\fg)$ (resp., 
$\Gr_-(\cF(\fg)_{int})$) with the subalgebra spanned by $a_{\nu}$ with $\nu\in P^+(\fg)$ (resp., 
$\nu\in P^+(\fg)_{int}$).

One has $a_0=1$; we denote by $\ev_0:\Gr_-(\CC)\to \mathbb{Z}$  the 
co-unit map given by $\ev_0(a_{\nu})=\delta_{0,\nu}$ for 
$\nu\in P^+(\fg)$. 
One has
\[
\Gr_-(\CC)=\Gr_-(\CC_{int})^W\oplus \Gr_-(\CC_{nint})^W,
\]
where $\CC_{nint}$ consists of modules with weights that lie in $P(\fg)$ and are non-integral. 
For $b\in \Gr_-(\CC)^W$ and $b'\in  \Gr_-(\CC_{nint})^W$ one has $bb'=\ev_0(b)b'$, so $ \Gr_-(\CC_{nint})$ is an ideal of $\Gr_-(\CC)$. 
By Corollary \ref{cor-ds-kernel}, for $x\not=0$, $\ds_x( \Gr_-(\CC_{nint})^W)=0$  and the image of $\ds_x$ lies in $\Gr_-(\CC_{int})$.
This reduces a study of $\Gr_-(\CC)$ to a study of $\Gr_-(\CC_{int})$.

\subsubsection{}
The ring  $\Gr_-(\cF(\fg)_{int})\otimes_{\mathbb{Z}} \mathbb{Z}[\sqrt{-1}]$ can be realized in the following way.

Let $V$ be a free $\mathbb{Z}[\sqrt{-1}]$-module with
a basis $\{v_i|\ i\in\mathbb{Z}\setminus\{0\}\}$.
Denote by $\bigwedge V$ the external ring of $V$. This is a $\mathbb{N}$-graded 
supercommutative ring; we consider $(\bigwedge V)^{spoil}$, which has
\[
\left(\bigwedge V\right)^{spoil}=\bigoplus\limits_{i=0}^{\infty} \bigwedge^{2i} V\oplus 
\bigoplus\limits_{i=0}^{\infty} \bigwedge^{2i+1} V\vareps
\]
where, we recall, $\vareps$ is a formal variable satisfying $\vareps^2=2\vareps=0$.
This is an $\mathbb{N}$-graded commutative and supercommutative ring (which means that 
$(\bigwedge V)^{spoil}_i(\bigwedge V)^{spoil}_j=0$ if $i,j$ are odd).

We denote by $\Xi_{int}$ the set of finite multisets $\{a_i\}_{i=1}^s$
with  $a_i\in\mathbb{Z}\setminus\{0\}$ and $a_i+a_j\not=0$ for all $1\leq i,j\leq s$.
For $A\in\Xi_{int}$ we denote by
$(\bigwedge V)^{spoil}_A$ the span of 
the elements $v_{i_1}\wedge v_{i_2}\wedge \ldots\wedge v_{i_p}$ 
(resp., the elements $\vareps v_{i_1}\wedge v_{i_2}\wedge \ldots\wedge v_{i_p}$) 
with
$\Core(\{i_1,\ldots,i_p\})=A$ if $A$ has an even (resp., an odd) cardinality.
Clearly, this gives a grading
\[
\left(\bigwedge V\right)^{spoil}_A\left(\bigwedge V\right)^{spoil}_B\sub\left(\bigwedge V\right)^{spoil}_{A\diamond B}
\]
and  $(\bigwedge V)^{spoil}_{\emptyset}$ is the subring of $(\bigwedge V)^{spoil}$ generated by $v_{i}\wedge v_{-i}$.   In other words $\left(\bigwedge V\right)^{spoil}$ is a $\Xi_{int}$-graded algebra.  For each $k$ we consider the ideal
$$J_k:=\displaystyle\sum_{i=k+1}^{\infty} \left(\bigwedge V\right)^{spoil}_i.$$

Recall that 
the map $\lambda\mapsto \Nonzero(\lambda)$ gives 
a one-to-one correspondence between $P^+(\fg)_{int}$ and the 
subsets $S\subset \mathbb{Z}\setminus\{0\}$ of cardinality at most $n$.
The weights in $I_0$ correspond to the subsets of 
even cardinality; the zero weight corresponds to $\emptyset$.
We define 
\[
\psi: \Gr_-(\cF(\fg)_{int})\otimes_{\mathbb{Z}}\mathbb{Z}[\sqrt{-1}]\iso
(\bigwedge V)^{spoil}/J_n
\]
by setting, for $\lambda\in I_i\cap P^+(\fg)$,
\[
\psi(a_{\lambda}):=\vareps^i\frac{t(\lambda)}{|t(\lambda)|}\cdot v_{\lambda},
\]
where $v_{\lambda}:=v_{i_1}\wedge v_{i_2}\wedge \ldots\wedge v_{i_k}$, and $\Nonzero(\lambda)=\{i_1>i_2>\ldots>i_k\}$.

Combining Corollary \ref{corollary-red-Gr-q(n)}, Proposition \ref{proposition-red-Gr-embedding-spoiled-polys} and Corollary \ref{cor-ds-q(n)} obtain

\begin{theorem}\label{thm 9.3}
\begin{enumerate}
\item  The map $\psi: \Gr_-(\cF(\fg)_{int})\otimes_{\mathbb{Z}}\mathbb{Z}[\sqrt{-1}]\iso \left(\bigwedge V\right)^{spoil}/J_n$ is an isomorphism of $\Xi_{int}$-graded spoiled super rings.

\item  In particular, for $A\in\Xi_{int}$, $\psi$ restricts to an isomorphism of vector spaces \newline $\psi_A:\Gr_-(\cF(\fg)_{A})\to\left(\bigwedge V\right)^{spoil}_A/J_n$.

\item The map $\ds_s: R(\fq_n)_{int}\to R(\fq_{n-s})_{int}$ corresponds to the natural quotient map 
\[
\left(\bigwedge V\right)^{spoil}/J_n\to \left(\bigwedge V\right)^{spoil}/J_{n-s}
\]
(one has $J_{n-s}\supset J_n$ for $s>0$).
\end{enumerate}
\end{theorem}

\begin{proof}
	It is clear that $\psi$ is an isomorphism of $\mathbb{Z}$-modules which preserves the $\Xi$-grading.  The result on $ds_s$ is also clear from the map.  Thus it remains to check that $\psi$ respects multiplication.    
	
	Let $\lambda,\mu\in P^+(\fg)$ with $\lambda,\mu\in I_0$; the case when either is in $I_1$ is much easier.  If $a_{\lambda}a_{\mu}=0$ then the result is clear, thus we assume that $a_{\lambda}a_{\mu}\neq0$.  In this case, there exists unique elements $\lambda'\in W\lambda$ and $\mu'\in W\mu$ such that $\lambda'+\mu'$ is dominant, and thus we have $a_{\lambda}a_{\mu}=(-1)^ia_{\lambda'+\mu'}$, where $E_{\lambda'}\otimes E_{\mu'}=\Pi^iE_{\lambda'+\mu'}$.  We may write $T_{\lambda'+\mu'}=(-1)^jT_{\lambda'}T_{\mu'}$ for some $j$, and it is clear that this $j$ also satisfies
	\[
	v_{\lambda'+\mu'}=(-1)^jv_{\lambda}v_{\mu}.
	\]
Then as in the proof of Proposition \ref{proposition-red-Gr-embedding-spoiled-polys}, we have $E_{\lambda'}\otimes E_{\mu'}=\Pi^{j+\ell}E_{\lambda'+\mu'}$, where $(-1)^\ell=\frac{t(\lambda)t(\mu)}{t(\lambda'+\mu')}\in\{\pm1\}$.  Therefore $a_{\lambda}a_{\mu}=(-1)^{j}\frac{t(\lambda)t(\mu)}{t(\lambda'+\mu')}a_{\lambda'+\mu'}$.  Therefore 
\[
\psi(a_{\lambda}a_{\mu})=(-1)^{j}\frac{t(\lambda)t(\mu)}{t(\lambda'+\mu'}\left(\frac{t(\lambda'+\mu')}{|t(\lambda'+\mu')|}v_{\lambda'+\mu'}\right)=\left(\frac{t(\lambda)}{|t(\lambda)|}v_{\lambda}\right)\left(\frac{t(\mu)}{|t(\mu)|}v_{\mu}\right)
\]
where we have used that $|t(\lambda'+\mu')|=|t(\lambda)||t(\mu)|$.
\end{proof}

\renewcommand{\arraystretch}{2}

For $A\in\Xi_{int}$, we define
\[
\ol{A}:=\left\{\begin{array}{ll}
\ \frac{\#\{i:a_i>0\}-\#\{i|a_i<0\}}{2} \text{ mod }2 & \text{ for } |A|\text{ even}; \\
\ \frac{\#\{i:a_i>0\}-\#\{i|a_i<0\}-1}{2} \text{ mod }2 & \text{ for } |A|\text{ odd}.
 \end{array} \right.
\]
Note that we have $\ol{A\diamond B}=\ol{A}+\ol{B}$ if $|A|$ or $|B|$ is even. 

\begin{lemma}
	The grading 
	\[
	\Gr_-(\cF(\fg)_{int})_{\ol{0}}=\bigoplus\limits_{\ol{A}=\ol{0}}\Gr_-(\cF(\fg)_A), \ \ \ \Gr_-(\cF(\fg)_{int})_{\ol{1}}=\bigoplus\limits_{\ol{A}=\ol{1}}\Gr_-(\cF(\fg)_A)
	\]
	defines another super ring structure on $\Gr_-(\cF(\fg)_{int})$; note that it does not have the structure of a spoiled super ring under this grading.
\end{lemma}

\begin{proof}
	This follows from the fact that $\Gr_-(\cF(\fg)_{int})$ is $\Xi_{int}$-graded, and $\Gr_-(\cF(\fg)_{A})\Gr_-(\cF(\fg)_{B})=0$ if $|A|$ and $|B|$ are odd.
\end{proof}
\begin{corollary}
	If $|A|$ is even, then $t(\lambda)\in\mathbb{R}^+$, and thus the isomorphism $\psi_A$ admits an integral structure.  In particular we have an isomorphism of graded rings:
\[
\bigoplus\limits_{\ol{A}=\ol{0}}\Gr_-(\cF(\fg)_A)\to \bigoplus\limits_{\ol{A}=\ol{0}}\left(\bigwedge V^{\mathbb{Z}}\right)_{A}^{spoiled}.
\]
Here $V^{\mathbb{Z}}$ denotes the free $\mathbb{Z}$-module with basis $\{v_i:i\in\mathbb{Z}\setminus\{0\}\}$.
\end{corollary}

\end{document}